\documentclass[12pt,letterpaper,fleqn]{article}

\title{On designing a resilient green supply chain to mitigate ripple effect: a two-stage stochastic optimization model}

\author{
 Hossein Mirzaee$^{1,2}$\\
  \texttt{hossein.mirzaee@usask.ca}
  \and
  Hamed Samarghandi$^{2}$\\
  \texttt{samarghandi@edwards.usask.ca}
  \and
  Keith Willoughby$^{2}$\\
  \texttt{willoughby@edwards.usask.ca}\blfootnote{1 Corresponding author}\blfootnote{2 Edwards School of Business, University of Saskatchewan, Saskatoon, SK, Canada} 
}

\usepackage{lipsum}

\newcommand\blfootnote[1]{%
  \begingroup
  \renewcommand\thefootnote{}\footnote{#1}%
  \addtocounter{footnote}{-1}%
  \endgroup
}

\usepackage{natbib}

\usepackage{soul}
\usepackage{attachfile}
\sethlcolor{lightgray}
\usepackage{cancel}
\usepackage[toc,page,title]{appendix}
\usepackage{adjustbox}
\usepackage[dvipsnames]{xcolor}
\usepackage{comment}
\usepackage{hyperref}
\hypersetup{colorlinks,linkcolor={blue},citecolor={blue},urlcolor={blue}}  
\usepackage{lscape}
\usepackage{amsmath}
\usepackage[left=1in,right=1in,top=1in,bottom=1in]{geometry}
\usepackage{xcolor}
\usepackage{booktabs}
\usepackage{float}
\usepackage{tikz}
\usepackage{pgfplotstable}
\usepackage[normalem]{ulem}
\usepackage{pgfplots}
\usepackage{url}
\usepackage{graphicx}
\usepackage[justification=centering]{caption}
\usetikzlibrary{positioning,calc}       
 \tikzstyle{block} = [draw,rectangle,thick,minimum height=2em,minimum width=2em,fill=green!10!white]
 \def\checkmark{\tikz\fill[scale=0.4](0,.35) -- (.25,0) -- (1,.7) -- (.25,.15) -- cycle;} 
\usepackage{graphicx}
\usepackage[english]{babel}
\date{}

\begin{document}

\maketitle
\begin{abstract}
\noindent Disasters and disruptions such as the COVID-19 pandemic can significantly interrupt supply chains and industries. To control these disruptions, decision-makers must focus on supply chain resiliency. This paper proposes a multi-stage, multi-period green supply chain design model and six resilience strategies, with downstream and upstream disruptions taken into account to analyze both the ripple and bullwhip effect, respectively. To control the mentioned disruptions and handle the uncertainties of parameter estimations, a two-stage stochastic optimization approach is devised. The objectives are to minimize the total cost of disruption, and $CO_{2}$ emission under the cap-and-trade mechanism as a government-issued emission regulation. The proposed decision-making framework and solution approach are validated using a numerical experiment followed by sensitivity analysis. The results show the optimum structure of the supply chain and the best resilient strategies to mitigate the ripple effect. Moreover, the effect of a decline in capacity of facilities on the optimal solution and the applied resilient strategies is investigated. This study provides managerial insights to help governments set the proper amount of cap, and supply chain managers to predict the demand behaviour of essential and non-essential products in the event of disruptions.

\vspace{1cm}
\end{abstract}

\setlength{\columnsep}{0.6cm}

\section{Introduction}\label{Introduction}

Disruptions caused by natural or human-made disasters affect supply chains in different aspects including transportation delays, labor unavailability, and supply-side shortage. A supply chain disruption announcement decreases a firm’s stock returns by 20\% on average after six months \citep{hendricks2005empirical}. Various examples demonstrate the challenges the firms face when trying to recover from a disruption: six months after Japan's tsunami in 2011, Toyota faced disruption in its supply network, and due to a shortage of parts, idled some of its plants in North America \citep{kim2015supply}. More recently, the COVID-19 pandemic outbreak caused long-term negative impacts on supply chains and revealed their vulnerabilities \citep{liu2022optimization}. These examples showcase the importance of adaptability and resiliency of supply chains in surviving new conditions in case of a sizeable disruption, which has recently gained attention among scholars and practitioners \citep{ivanov2022stress}.

One type of interruption to scrutinize for improving supply chain adaptability is the ripple effect, which is described as the propagation of disturbances that arise from the disruption of supply chain elements \citep{ivanov2016disruption}. The adverse impacts of the ripple effect spread downstream in the supply chain \citep{monostori2021mitigation}. Real-world examples emphasize that controlling the ripple effect is crucial for supply chain managers. For instance, in June 2020, Mercedes-Benz ceased production of an off-road vehicle in Alabama as a result of a shortage in components imported from its European suppliers during the global COVID-19 pandemic \citep{REUTERS2020news}. 

The desirable approach for efficient recovery from the impact of ripple effect is constructing intrinsic supply chain resiliency. Having contingency plans such as backup suppliers or temporary facilities at the supply chain design stage is helpful in controlling the ripple effect \citep{ivanov2015supply}. In other words, appropriate strategies must be considered during the design stage to mitigate the crunch in the aftermath of inadmissible events such as supply delay, demand hike, or capacity contraction \citep{sharma2022managing}. The auspicious design strategies include, but are not limited to, considering backup suppliers, capacity expansion and multiple assignments \citep{gholami2021design}.

Another important factor in designing a supply chain is attention to the environmental aspects as they bring competitive advantages for the firms \citep{boskabadi2022design}. Devising emission abatement schemes, producing recyclable products, and using green technologies are some of the elements that lead to a greener supply chain \citep{mirzaee2022three}.

On another note, inaccurate estimation of the design parameters may result in colossal losses in uncertain environments \citep{wang2021incentive}. Ergo, uncertainty in the forecast values of the parameters is another factor that negatively impacts the supply chain performance, which necessitates adopting an appropriate approach. The three most common uncertainty control methods are stochastic optimization, robust optimization, and fuzzy optimization \citep{tordecilla2021simulation}, among which stochastic optimization is the most popular technique in the ripple effect literature. 

Stochastic optimization takes into account disruptions by using scenario-based modeling \citep{oksuz2020two}, while robust optimization, despite its several advantages, focuses directly on the worst-case scenario, which is not always relevant \citep{ivanov2019impact}. Conversely, fuzzy optimization prevents considering some scenarios regarding the ripple effect \citep{ozccelik2021robust}. Moreover, it requires deep knowledge about the problem's parameters to develop a membership function, which is not always applicable \citep{memon2015group}. Henceforth, to provide the best strategic and operational decisions, this paper adopts a two-stage stochastic optimization approach to control uncertainty. 

This study aims to address the following research questions:
\begin{itemize}
    \item What are the strategic and operational decisions to achieve a green, sustainable and resilient supply chain?
    \item What are the strategies that best mitigate the ripple effect?
    \item How can the decision makers best control the inherent uncertainty which is innate to the estimation of the parameter values? 
\end{itemize}

The remainder of this paper is structured as follows: a summary of the relevant research is provided in section \ref{Literature review}. Section \ref{Model formulation} explains the developed two-stage stochastic optimization model for the resilient, green supply chain design problem. Section \ref{Numerical experiment} is devoted to numerical experiments and their analysis. Section \ref{Conclusion} presents the concluding remarks and the venues for future research.

\section{Literature review}\label{Literature review}

The literature related to this study focuses on three aspects. First is the environmental perspectives of the supply chain design. The second facet is linked to resilience strategies of mitigating the ripple effect. The third prospect is associated with employing stochastic optimization to deal with parameter uncertainty and ripple effect. Table \ref{Table: LR} lists the related studies and highlights their opposing views. The rest of this section recaps the corresponding literature, the existing gaps, and the contributions of this paper.

Green supply chain design integrates environmental issues into strategic decisions \citep{foroozesh2022green}. The significance of these issues has attracted the attention of various scholars recently \citep{bhatia2021green}. The role of environmental investment in the supply chain network configuration phase in making the supply chain greener is the main focus of \cite{wang2011multi}. In their study, they considered $CO_{2}$ emission as the main indicator of supply chain greenness, which is a mainstream index for environmental issues and can be estimated with more ease. In a similar study \cite{o2013fifty} considered environmental issues in GSCD by explaining limited goals such as GHG emission. \cite{hasani2021multi} take environmental and economic concerns into account in designing a supply chain network. They utilized resilient strategies to mitigate disruptions. \cite{mohebalizadehgashti2020designing} formulated a multi-objective MILP for the GSCD problem which aims to minimize $CO_{2}$ emissions from transportation and maximize total utilization of facility capacities. 

In recent years, especially after the COVID-19 pandemic, supply chain resiliency has received extensive attention. For instance, several researchers published review papers to comprehensively study resilience strategies in supply chains \citep{hosseini2019review, ivanov2019impact, snyder2016or}. \cite{tomlin2006value} specify that resilience strategies are classified into two main groups: pre-disruption and post-disruption schemes. For instance, accumulating safety stock is a pre-disruption resilience strategy for a situation when the supply side of the supply chain is affected \citep{foroozesh2022green}. \cite{yilmaz2021ensuring} introduced four stages of controlling the ripple effect, namely, preparation, first response, preparation for recovery, and recovery. They suggest employing pre-disruption resilience strategies in the first three stages and utilizing a post-disruption resilient strategy for the last stage. 

\cite{ni2018modeling} recommend applying post-disruption strategies such as using contingency supplies furnished by backup suppliers or stockpiling systems to maintain customer satisfaction and responding to unmet demand. \cite{kamalahmadi2016developing} develop a two-stage MIP to design a sourcing plan with high flexibility. They combine the transportation channel selection problem and supplier selection and order allocation problem, and devise contingency plans for mitigating the negative effects of disruption to minimize total supply chain cost. They found that contingency plan implementation increases suppliers' flexibility in adapting to manufacturers' capacity, and reduces disruption's severity. \cite{jabbarzadeh2018resilient} applied and assessed inventory levels and backup suppliers as two resilient strategies. \cite{hosseini2020ripple} considered segregating suppliers as a resilience strategy in supply chain design.

One way to prevent the spread of the ripple effect and to disallow parameter uncertainty to negatively affect the predictions is utilizing a proper uncertainty control method. Accoding to \cite{rezapour2017resilient} and \cite{kamalahmadi2016developing}, robust, stochastic, and fuzzy techniques are more prominent than the deterministic models to cope with uncertainties. \cite{badri2017two} developed a two-stage stochastic optimization model to maximize the total value of a supply chain. \cite{yilmaz2021ensuring} applied a two-stage stochastic technique to design a reverse supply chain in the presence of ripple effect, and showed that, as a result, the emission level increases by 40\%. Therefore, emission abatement regulations should be enforced to avoid the upsurge. In this paper, a two-stage stochastic optimization approach is utilized to cope with uncertainties. For more information about multi-stage stochastic optimization and implementation of stochastic techniques, the interested reader is referred to \cite{khaloie2020coordinated} and \cite{cui2020stochastic}.

Table \ref{Table: LR} summarizes the recent and relevant studies. According to the literature summary and the above-mentioned papers, one notices that although there is more emphasis on minimizing the GHG emissions in green supply chain management, environmental considerations are not the main focus of supply chain design problems, and governmental emission reduction regulations such as cap-and-trade are considered as a hindrance that work against maximizing the profits. Furthermore, the ripple effect, as one of the main disruption elements, is rarely considered in resilient supply chain design studies. 

The papers that have studied the ripple effect have mainly considered resilience strategies in their general form, such as pre-disruption and post-disruption schemes; specific resilience strategies such as safety stock and backup suppliers have rarely been studied. In other words, resilience strategies are more studied based on pre- or post-disruption classification, and the details of particular strategies like temporary facilities and safety stock are neglected. Regarding the uncertainty control approaches, most of the papers have used robust (RO) or fuzzy optimization (FO) and stochastic programming (SP); few studies consider multi-stage stochastic programming to handle uncertainty. The literature related to this study is limited to the general form of disruption and rarely studies the areas of the supply chain affected by disruptions. To the best of our knowledge, there is no study that considers both upstream and downstream propagation of disruptions in a supply chain.

This paper studies the ways to mitigate the ripple effect and demand uncertainty by developing a multi-period, multi-stage green resilient supply chain. We consider six resilience strategies to keep the supply and demand side of the supply chain in control. We deploy a two-stage stochastic optimization approach as an effective way of controlling parameter estimation uncertainty and ripple effect.

\begin{table}
\centering
\caption{Literature review summary}
\begin{adjustbox}{max width=\textwidth}
\renewcommand{\arraystretch}{1}
\newcommand*{\TitleParbox}[1]{\parbox[c]{10.75cm}{\raggedright #1}}%
\begin{tabular}{|lccccccc|}
\hline
\multicolumn{1}{|l|}{\textbf{Authors}}                                                  & \multicolumn{1}{l|}{\begin{tabular}[c]{@{}c@{}}\textbf{Network} \\   \textbf{stages}\end{tabular}} & \multicolumn{1}{l|}{\begin{tabular}[c]{@{}l@{}}\textbf{Multi-}\\ \textbf{period}\end{tabular}} & \multicolumn{1}{c|}{\begin{tabular}[c]{@{}c@{}}\textbf{Green} \\ \textbf{supply chain}\end{tabular}} & \multicolumn{1}{c|}{\begin{tabular}[c]{@{}c@{}}\textbf{The ripple} \\      \textbf{effect}\end{tabular}} & \multicolumn{1}{c|}{\textbf{Uncertainty}} & \multicolumn{1}{c|}{\textbf{Resilient strategies}}                                                                                                                                     & \textbf{Disruption effects}                  \\ \hline
\multicolumn{1}{|l|}{\cite{kamalahmadi2016developing}}       & \multicolumn{1}{c|}{2}                                                           & \multicolumn{1}{c|}{\_}                                                      & \multicolumn{1}{c|}{\_}                                                            & \multicolumn{1}{c|}{\_}                                                                & \multicolumn{1}{c|}{Others}              & \multicolumn{1}{c|}{Contingency plans}                                                                                                                                        & Supply disruption                   \\ \hline
\multicolumn{1}{|l|}{\cite{fattahi2017responsive}}     & \multicolumn{1}{c|}{3}                                                           & \multicolumn{1}{c|}{\checkmark}                               & \multicolumn{1}{c|}{\_}                                                            & \multicolumn{1}{c|}{\_}                                                                & \multicolumn{1}{c|}{Two-stage SP}        & \multicolumn{1}{c|}{\begin{tabular}[c]{@{}c@{}}Contingency and \\ mitigation(general)\end{tabular}}                                                                           & Facility capacity                \\ \hline
\multicolumn{1}{|l|}{\cite{mohammed2017cost}}                & \multicolumn{1}{c|}{3}                                                           & \multicolumn{1}{c|}{\_}                                                      & \multicolumn{1}{c|}{\_}                                                            & \multicolumn{1}{c|}{\_}                                                                & \multicolumn{1}{c|}{SP}                  & \multicolumn{1}{c|}{\_}                                                                                                                                                       & \_                                  \\ \hline
\multicolumn{1}{|l|}{\cite{rezapour2017resilient}}          & \multicolumn{1}{c|}{3}                                                           & \multicolumn{1}{c|}{\_}                                                      & \multicolumn{1}{c|}{\_}                                                            & \multicolumn{1}{c|}{\checkmark}                                         & \multicolumn{1}{c|}{Deterministic}       & \multicolumn{1}{c|}{\begin{tabular}[c]{@{}c@{}}Emergency stock at retailers,   \\ backup capacities at \\ suppliers, multiple sourcing\end{tabular}}                          & Supplier disruption               \\ \hline
\multicolumn{1}{|l|}{\cite{pavlov2017hybrid}}                & \multicolumn{1}{c|}{4}                                                           & \multicolumn{1}{c|}{\_}                                                      & \multicolumn{1}{c|}{\_}                                                            & \multicolumn{1}{c|}{\checkmark}                                         & \multicolumn{1}{c|}{FO}                  & \multicolumn{1}{c|}{\_}                                                                                                                                                       & \_                                  \\ \hline
\multicolumn{1}{|l|}{\cite{badri2017two}}                     & \multicolumn{1}{c|}{3}                                                           & \multicolumn{1}{c|}{\checkmark}                               & \multicolumn{1}{c|}{\_}                                                            & \multicolumn{1}{c|}{}                                                                  & \multicolumn{1}{l|}{Two-stage SP}        & \multicolumn{1}{c|}{\_}                                                                                                                                                       & \_                                  \\ \hline
\multicolumn{1}{|l|}{\cite{amiri2018multi}}     & \multicolumn{1}{c|}{3}                                                           & \multicolumn{1}{c|}{\checkmark}                               & \multicolumn{1}{c|}{\_}                                                            & \multicolumn{1}{c|}{\_}                                                                & \multicolumn{1}{c|}{Two-stage SP}        & \multicolumn{1}{c|}{\_}                                                                                                                                                       & \multicolumn{1}{c|}{\_}               \\ \hline
\multicolumn{1}{|l|}{\cite{zahiri2018design}}          & \multicolumn{1}{c|}{3}                                                           & \multicolumn{1}{c|}{\checkmark}                               & \multicolumn{1}{c|}{\_}                                                            & \multicolumn{1}{c|}{\_}                                                                & \multicolumn{1}{c|}{FO}                  & \multicolumn{1}{c|}{\_}                                                                                                                                                       & \multicolumn{1}{c|}{\_}               \\ \hline
\multicolumn{1}{|l|}{\cite{john2018multi}}                    & \multicolumn{1}{c|}{1}                                                           & \multicolumn{1}{c|}{\checkmark}                               & \multicolumn{1}{c|}{\_}                                                            & \multicolumn{1}{c|}{\_}                                                                & \multicolumn{1}{c|}{Deterministic}       & \multicolumn{1}{c|}{\_}                                                                                                                                                       & \_                                  \\ \hline
\multicolumn{1}{|l|}{\cite{liao2018reverse}}                 & \multicolumn{1}{c|}{2}                                                           & \multicolumn{1}{c|}{\_}                                                      & \multicolumn{1}{c|}{\checkmark}                                     & \multicolumn{1}{c|}{\_}                                                                & \multicolumn{1}{c|}{Deterministic}       & \multicolumn{1}{c|}{\_}                                                                                                                                                       & \_                                  \\ \hline
\multicolumn{1}{|l|}{\cite{ni2018modeling}}                   & \multicolumn{1}{c|}{1}                                                           & \multicolumn{1}{c|}{\checkmark}                               & \multicolumn{1}{c|}{\_}                                                            & \multicolumn{1}{c|}{\_}                                                                & \multicolumn{1}{l|}{Two-stage SP}        & \multicolumn{1}{c|}{\begin{tabular}[c]{@{}c@{}}Backup facilities, safety stock,   \\ idle capacity reserve\end{tabular}}                                                      & Demand                              \\ \hline
\multicolumn{1}{|l|}{\cite{jabbarzadeh2018resilient}}         & \multicolumn{1}{c|}{2}                                                           & \multicolumn{1}{c|}{\_}                                                      & \multicolumn{1}{c|}{\checkmark}                                     & \multicolumn{1}{c|}{\_}                                                                & \multicolumn{1}{c|}{SP}                  & \multicolumn{1}{c|}{\begin{tabular}[c]{@{}c@{}}Backup suppliers, production   \\ capacity expansion\end{tabular}}                                                             & Supply disruption                   \\ \hline
\multicolumn{1}{|l|}{\cite{sawik2019disruption}}        & \multicolumn{1}{c|}{2}                                                           & \multicolumn{1}{c|}{\checkmark}                               & \multicolumn{1}{c|}{\_}                                                            & \multicolumn{1}{c|}{\_}                                                                & \multicolumn{1}{c|}{Two-stage SP}   & \multicolumn{1}{c|}{\_}                                                                                                                                                       & \multicolumn{1}{c|}{\_}               \\ \hline
\multicolumn{1}{|l|}{\cite{darestani2019robust}}              & \multicolumn{1}{c|}{3}                                                           & \multicolumn{1}{c|}{\checkmark}                               & \multicolumn{1}{c|}{\checkmark}                                     & \multicolumn{1}{c|}{\_}                                                                & \multicolumn{1}{c|}{RO}                  & \multicolumn{1}{c|}{\_}                                                                                                                                                       & \_                                  \\ \hline
\multicolumn{1}{|l|}{\cite{hosseini2019new}}                 & \multicolumn{1}{c|}{2}                                                           & \multicolumn{1}{c|}{\_}                                                      & \multicolumn{1}{c|}{\_}                                                            & \multicolumn{1}{c|}{\checkmark}                                         & \multicolumn{1}{c|}{Others}              & \multicolumn{1}{c|}{\_}                                                                                                                                                       & \_                                  \\ \hline
\multicolumn{1}{|l|}{\cite{hosseini2019reverse}}             & \multicolumn{1}{c|}{3}                                                           & \multicolumn{1}{c|}{\_}                                                      & \multicolumn{1}{c|}{\checkmark}                                     & \multicolumn{1}{c|}{\_}                                                                & \multicolumn{1}{c|}{Others}              & \multicolumn{1}{c|}{\_}                                                                                                                                                       & \_                                  \\ \hline
\multicolumn{1}{|l|}{\cite{zhang2019fuzzy}}                  & \multicolumn{1}{c|}{2}                                                           & \multicolumn{1}{c|}{\_}                                                      & \multicolumn{1}{c|}{\_}                                                            & \multicolumn{1}{c|}{\checkmark}                                         & \multicolumn{1}{c|}{FO}                  & \multicolumn{1}{c|}{\begin{tabular}[c]{@{}c@{}}Backup manufacturer, multiple   \\ distributor\end{tabular}}                                                                   & Supply disruption                   \\ \hline
\multicolumn{1}{|l|}{\cite{hosseini2020ripple}}               & \multicolumn{1}{c|}{2}                                                           & \multicolumn{1}{c|}{\_}                                                      & \multicolumn{1}{c|}{\_}                                                            & \multicolumn{1}{c|}{\checkmark}                                         & \multicolumn{1}{c|}{SP}                  & \multicolumn{1}{c|}{\_}                                                                                                                                                       & Supplier disruption               \\ \hline
\multicolumn{1}{|l|}{\cite{tucker2020incentivizing}}         & \multicolumn{1}{c|}{3}                                                           & \multicolumn{1}{c|}{\checkmark}                               & \multicolumn{1}{c|}{\_}                                                            & \multicolumn{1}{c|}{\checkmark}                                         & \multicolumn{1}{c|}{SP}                  & \multicolumn{1}{c|}{\begin{tabular}[c]{@{}c@{}}Configuration of suppliers and   \\ manufacturers, safety stock\end{tabular}}                                                  & Supply disruption                   \\ \hline
\multicolumn{1}{|l|}{\cite{mohebalizadehgashti2020designing}} & \multicolumn{1}{c|}{3}                                                           & \multicolumn{1}{c|}{\checkmark}                               & \multicolumn{1}{c|}{\checkmark}                                     & \multicolumn{1}{c|}{\_}                                                                & \multicolumn{1}{c|}{Others}              & \multicolumn{1}{c|}{\_}                                                                                                                                                       & \_                                  \\ \hline
\multicolumn{1}{|l|}{\cite{ozccelik2021robust}}             & \multicolumn{1}{c|}{2}                                                           & \multicolumn{1}{c|}{\_}                                                      & \multicolumn{1}{c|}{\checkmark}                                     & \multicolumn{1}{c|}{\checkmark}                                         & \multicolumn{1}{c|}{RO}                  & \multicolumn{1}{c|}{\_}                                                                                                                                                       & \_                                  \\ \hline
\multicolumn{1}{|l|}{\cite{hasani2021multi}}                 & \multicolumn{1}{c|}{3}                                                           & \multicolumn{1}{c|}{\checkmark}                               & \multicolumn{1}{c|}{\checkmark}                                     & \multicolumn{1}{c|}{\_}                                                                & \multicolumn{1}{c|}{RO}                  & \multicolumn{1}{c|}{\begin{tabular}[c]{@{}c@{}}Backup suppliers, facility   \\ dispersion, facility fortification\end{tabular}}                                               & Supply disruption                   \\ \hline
\multicolumn{1}{|l|}{\cite{yilmaz2021ensuring}}              & \multicolumn{1}{c|}{3}                                                           & \multicolumn{1}{c|}{\_}                                                      & \multicolumn{1}{c|}{\checkmark}                                     & \multicolumn{1}{c|}{\checkmark}                                         & \multicolumn{1}{l|}{Two-stage SP}        & \multicolumn{1}{c|}{Temporary facilities}                                                                                                                                     & Supply disruption                   \\ \hline
\multicolumn{1}{|l|}{\cite{foroozesh2022green}}             & \multicolumn{1}{c|}{4}                                                           & \multicolumn{1}{c|}{\checkmark}                               & \multicolumn{1}{c|}{\checkmark}                                     & \multicolumn{1}{c|}{\_}                                                                & \multicolumn{1}{c|}{FO}                  & \multicolumn{1}{c|}{\begin{tabular}[c]{@{}c@{}}Multiple sourcing, \\ horizontal collaboration,   \\ coverage radius\end{tabular}}                                             & Supply disruption                   \\ \hline
\multicolumn{1}{|l|}{This study}                                               & \multicolumn{1}{c|}{4}                                                           & \multicolumn{1}{c|}{\checkmark}                               & \multicolumn{1}{c|}{\checkmark}                                     & \multicolumn{1}{c|}{\checkmark}                                         & \multicolumn{1}{l|}{Two-stage SP}        & \multicolumn{1}{c|}{\begin{tabular}[c]{@{}c@{}}Backup suppliers, multiple   \\ sourcing, temporary \\ facilities, blockchain, \\      safety stock, stockpiling\end{tabular}} & \multicolumn{1}{c|}{Supply, demand} \\ \hline
\multicolumn{8}{|l|}{SP=Stochastic   programming, RO= Robust optimization, FO=Fuzzy optimization}                                                                                \\ \hline
\end{tabular}
\end{adjustbox}
\label{Table: LR}
\end{table}

\section{Resilient Green Supply Chain Design (RGSCD) model formulation}\label{Model formulation}

\subsection{Problem statement}\label{Problem statement}

\begin{figure}[h]
\centering
\includegraphics[width=18cm]{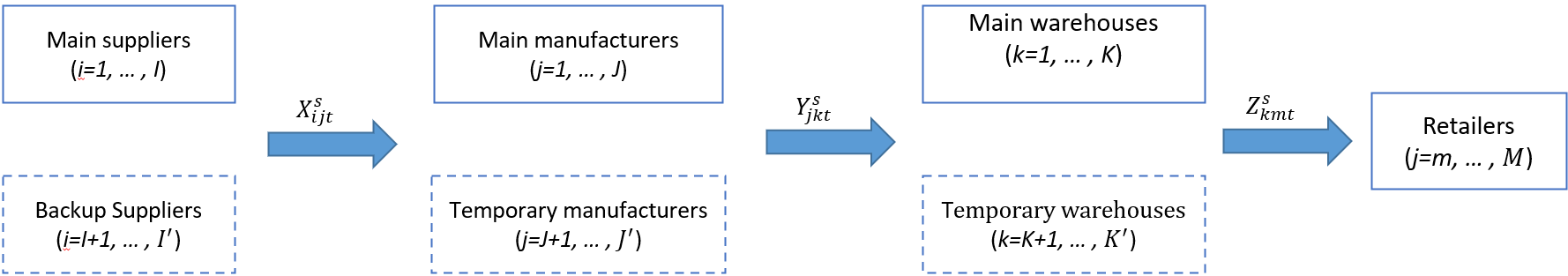}
\caption{Supply chain structure} 
\label{Fig: SC structure}
\end{figure}

The resilient green supply chain design problem can be stated as follows: there are four stages represented by sets of $I$, $J$, $K$, and $M$, i.e., suppliers, manufacturers, warehouses, and retailers, respectively. Raw material is procured to manufacturers by suppliers ($X_{ijt}^{s}$). After producing products, they are kept in warehouses ($T_{jkt}^{s}$) to be sent to retailers ($Z_{kmt}^{s}$). There are $L$ different transportation modes to move the materials and products between the supply chain echelons. $T$ time periods are defined for this problem to make sure the right decisions are made considering long-term planning period. Different levels of the ripple effect are shown by $S$ possible scenarios. Backup suppliers, and temporary manufacturing and warehousing facilities are considered as extra sources in the event of disruption. In order to satisfy a predetermined service level, these backup suppliers and facilities are reserved to be utilized when the ripple effect emerges in the supply chain. The supply chain structure is depicted in figure \ref{Fig: SC structure}. The objective is to minimize the costs of supply chain operational decisions, strategic resilient decisions, and total $CO_{2}$ emission during production and transportation process of a product. It is assumed that a cap-and-trade mechanism is in place as an emission abatement mechanism.

The developed model in this paper is formulated based on the following assumptions:

\begin{itemize} 
    \item $CO_{2}$ is emitted as a result of transportation as well as production processes in the supply chain.
    \item The amount of $CO_{2}$ emission is restricted by government.
    \item Different possible scenarios are considered for the uncertain parameters to show the disparate levels of disruptions.
    \item The manufacturers are allowed to have shortage, which will be back-ordered.
    \item All temporary manufacturing centers and warehouses are identical in terms of production and storage capacity.
\end{itemize} 

\subsection{Two-stage stochastic optimization model}\label{Model:Two-stage stochastic optimization model}

This section presents a two-stage stochastic mixed integer linear optimization model with the aim of minimizing total costs and emissions. The main goal of the model is to develop a resilient supply chain that is capable of controlling the ripple effect. The proposed model ensures that emission levels are controlled such that the governmental issued emission restrictions are met.

The proposed two-stage stochastic optimization model is as follows.

\textbf{Indices and sets}

\noindent $i$: Main and backup suppliers $(i=1,\cdots, I, I+1, \cdots I^{\prime})$

\noindent $j$: Permanent and temporary manufacturing facilities $(j=1,\cdots, J, J+1, \cdots J^{\prime})$

\noindent $k$: Permanent and temporary warehouses $(k=1,\cdots, K, K+1, \cdots K^{\prime})$

\noindent $m$: Retailers $(m=1, \cdots, M)$

\noindent $t$: Time periods $(t=1,\cdots,T)$

\noindent $l$: Transportation modes $(l=1, \cdots, L)$

\noindent $s$: Scenarios $(s=1,\cdots,S)$

\textbf{Parameters}

\noindent $TC1_{ijl}^{s}$: Cost of transportation from supplier $i$ to manufacturer $j$ using transportation mode $l$ under scenario $s$

\noindent $TC2_{jkl}^{s}$: Cost of transportation from manufacturer $j$ to warehouse $k$ using transportation mode $l$ under scenario $s$

\noindent $TC3_{kml}^{s}$: Cost of transportation from warehouse $k$ to retailer $m$ using transportation mode $l$ under scenario $s$

\noindent $T1_{ijl}^{s}\ge0$: Transportation delay from supplier $i$ to manufacturer $j$ using transportation mode $l$ for one order under scenario $s$

\noindent $T2_{jkl}^{s}\ge0$: Transportation delay from manufacturer $j$ to warehouse $k$ using transportation mode $l$ for one order under scenario $s$

\noindent $T3_{kml}^{s}\ge0$: Transportation delay from warehouse $k$ to retailer $m$ using transportation mode $l$ for one order under scenario $s$

\noindent $Ms_{j}$: Setup cost of opening temporary manufacturing center $j \in \{j=J+1, \cdots, J^{\prime}\}$

\noindent $Ws_{k}$: Setup cost of opening temporary warehouse $k\in \{K+1, \cdots, K^{\prime}\}$ 

\noindent $So$: Setup cost of establishing an information sharing system

\noindent $St1_{ijl}$: Saved ordering time due to using the information sharing system for an order placed by manufacturer $j$ to supplier $i$, per order

\noindent $St2_{jkl}$: Saved ordering time due to using the information sharing system for an order placed by warehouse $k$ to manufacturer $j$, per order

\noindent $St3_{kml}$: Saved ordering time due to using the information sharing system for an order placed by retailer $m$ to warehouse $k$, per order

\noindent $Dc$: Delay cost per order per day

\noindent $Tr$: Training cost of employees to operate the information sharing system

\noindent $Mic_{j}^{s}$: Inventory cost per product for manufacturer $j$ under scenario $s$

\noindent $Msc_{k}^{s}$: Shortage cost per product for manufacturer $j$ under scenario $s$

\noindent $D_{mt}^{s}$: Demand received by retailer $m$ in period $t$ under scenario $s$

\noindent $Mc_{j}$: Production capacity of manufacturer $j$

\noindent $Wc_{k}$: Capacity of warehouse $k$

\noindent $Spp$: Extra cost of stockpiled products

\noindent $Sm$: Minimum number of suppliers

\noindent $Cep_{j}$: Carbon emission of manufacturer $j$ during production process per product

\noindent $Cet1_{ijl}$: Carbon emission of transportation from  supplier $i$ to manufacturing center $j$ using transportation mode $l$, per product

\noindent $Cet2_{jkl}$: Carbon emission of transportation from manufacturer $j$ to warehouse $k$ using transportation mode $l$, per product

\noindent $Cet3_{kml}$: Carbon emission of transportation from warehouse $k$ to retailer $m$ using transportation mode $l$, per product 

\noindent $Rce$: Reduced amount of carbon emission per unit of product shipped due to using information sharing system 

\noindent $Cap_{t}$: Maximum allowed carbon emission in period $t$

\noindent $r_{j}^{s}$: Reduced capacity ratio of manufacturer $j$ under scenario $s$ due to the ripple effect

\noindent $r^{'s}_{k}$: Reduced capacity ratio of warehouse $k$ under scenario $s$ due to the ripple effect

\noindent $M$: A large positive number

\noindent $Pr^{s}$: Occurrence probability of scenario $s$

\textbf{First stage decision variables}

\noindent $XX_{ijtl}$: A binary variable showing material flow between manufacturer $j$ and supplier $i$ in period $t$ using transportation mode $l$; $XX_{ijtl}=1$ means an order is received by manufacturer $j$ from supplier $i$

\noindent $YY_{jktl}$: A binary variable showing product flow between manufacturer $j$ and warehouse $k$ using transportation mode $l$ in period $t$; $YY_{jktl}=1$ means product flow between manufacturer and warehouse

\noindent $ZZ_{kmtl}$: A binary variable showing product flow between retailer $m$ and warehouse $k$ using transportation mode $l$ in period $t$; $ZZ_{kmtl}=1$ means product flow between warehouse and retailer

\textbf{Second stage decision variables}

\noindent $X_{ijtl}^{s}$: An integer variable showing the number of products transported from supplier $i$ to manufacturer $j$ in period $t$ using transportation mode $l$ under scenario $s$

\noindent $Y_{jktl}^{s}$: An integer variable showing the number of products transported from manufacturer $j$ to warehouse $k$ in period $t$ using transportation mode $l$ under scenario $s$

\noindent $Z_{kmtl}^{s}$: An integer variable showing the number of products transported from warehouse $k$ to retailer $m$ in period $t$ using transportation mode $l$ under scenario $s$

\noindent $MI_{jt}^{s}$: The $j^{th}$ manufacturer’s inventory at the end of period $t$ under scenario $s$

\noindent $MS_{jt}^{s}$: The $j^{th}$ manufacturer’s shortage at the end of period $t$ under scenario $s$

\noindent $MSS_{jt}^{s}$: The amount of $j^{th}$ manufacturer’s safety stock in period $t$ under scenario $s$

\noindent $SP_{jt}^{s}$: The number of ordered stockpiled products by manufacturer $j$ in period $t$ under scenario $s$

\textbf{Mathematical formulation}

This section presents a two-stage stochastic optimization model with two objectives. The first objective introduces all the supply chain costs; the second objective belongs to the $CO_{2}$ emissions. The first stage considers the deterministic variables, which are not dependent to the uncertain parameters. The second stage takes the uncertain variables into account. In this paper, uncertainty is handled by considering different possible scenarios that may occur due to disruption. The first stage optimizes the channel through which the raw materials transform to the  final product. In the second stage, $s$ scenarios are defined for the uncertain variables. The optimal solution of the presented model is obtained based on all possible scenarios in accordance with their occurrence probabilities.

The supply chain costs are as follows:
\begin{equation}\label{eq1: Transportation cost}
\begin{split}   
  C_{1} = \sum_{i,i^{'}\in I} \sum_{j,j^{'}\in J} \sum_{t\in T} \sum_{l\in L} \sum_{s\in S} XX_{ijtl}Dc(Tl1_{ijl}^{s}-St1_{ijl}) \\ +\sum_{j,j^{'}\in J} \sum_{k,k^{'}\in K} \sum_{t\in T} \sum_{l\in L} \sum_{s\in S} YY_{jktl}Dc(Tl2_{jkl}^{s}-St2_{jkl}) \\ + \sum_{k,k^{'}\in K} \sum_{m \in M} \sum_{t\in T} \sum_{l\in L} \sum_{s\in S} ZZ_{kmtl}Dc(Tl3_{kml}^{s}-St3_{kml})
\end{split}
\end{equation} 

Equation \eqref{eq1: Transportation cost} shows the cost of transportation delay between different supply chain echelons.
\begin{equation}\label{eq2: Transportation cost}
\begin{split}   
   C_{2}^{s} = \sum_{i,i^{'}\in I} \sum_{j,j^{'}\in J} \sum_{t\in T} \sum_{l\in L} \sum_{s\in S} X_{ijtl}^{s}Tc1_{ijl}^{s} + \sum_{j,j^{'}\in J} \sum_{k,k^{'}\in K} \sum_{t\in T} \sum_{l\in L} \sum_{s\in S} Y_{jktl}^{s}Tc2_{jkl}^{s} \\ + \sum_{k,k^{'}\in K} \sum_{m\in M} \sum_{t\in T} \sum_{l\in L} \sum_{s\in S} Z_{kmtl}^{s}Tc3_{kml}^{s}
\end{split}
\end{equation} 

Equation \eqref{eq2: Transportation cost} refers to the transportation cost throughout the supply chain. Equations \eqref{eq3: Inventory cost} and \eqref{eq4: Shortage cost} indicate the cost of holding inventory and the cost of shortage at the end of each period, respectively.
\begin{equation}\label{eq3: Inventory cost}
\begin{split}   
  C_{3}^{s} = \sum_{j,j^{'}\in J} \sum_{t\in T} \sum_{s\in S} MI_{jt}^{s}Mic_{j}^{s}
\end{split}
\end{equation} 
\begin{equation}\label{eq4: Shortage cost}
\begin{split}   
  C_{4}^{s} = \sum_{j,j^{'}\in J} \sum_{t\in T} \sum_{s\in S} MS_{jt}^{s}Msc_{j}^{s}
\end{split}
\end{equation}

The amount of $CO_{2}$ emission during the transportation and production processes is calculated by equation \eqref{eq5: Emission}.
\begin{equation}\label{eq5: Emission}
\begin{split}   
  E^{s} = \sum_{j,j^{'}\in J} \sum_{k,k^{'}\in K} \sum_{t\in T} \sum_{l\in L} \sum_{s\in S} Y_{jktl}^{s}Cep_{j} + \sum_{i,i^{'}\in I} \sum_{j,j^{'}\in J} \sum_{t\in T} \sum_{l\in L} \sum_{s\in S} X_{ijtl}^{s}Cet1_{ijl}\quad \quad \quad\quad\\ +  \sum_{j,j^{'}\in J} \sum_{k,k^{'}\in K} \sum_{t\in T} \sum_{l\in L} \sum_{s\in S} Y_{jktl}^{s}Cet2_{jkl} + \sum_{k,k^{'}\in K} \sum_{m,m^{'}\in M} \sum_{t\in T} \sum_{l\in L} \sum_{s\in S} Z_{kmtl}^{s}Cet3_{kml}
\end{split}
\end{equation} 

In the two-stage stochastic optimization model, the objective functions are formulated as follows:
\begin{equation}\label{eq6}
\begin{split}   
  Min\,Z_{1} = C_{1} + \sum_{s\in S}Pr^{s}(C_{2}^{s}+C_{3}^{s}+C_{4}^{s})
\end{split}
\end{equation} 
\begin{equation}\label{eq7}
\begin{split}   
  Min\,Z_{2} = \sum_{s\in S}Pr^{s}E^{s}
\end{split}
\end{equation} 

The first objective function \eqref{eq6} minimizes the total cost of supply chain; the second objective \eqref{eq7} minimizes total $CO_{2}$ emission. In these objective functions the optimal value of uncertain variables are obtained based on all possible scenarios.

\textbf{Constraints}

Constraints \eqref{eq8: }-\eqref{eq10: } specify which supply chain partners are used.
\begin{align}
    X_{ijtl}^{s} \leq XX_{ijtl}M \quad \qquad \qquad \forall  i,j,t,l,s  \label{eq8: }\\
    Y_{jktl}^{s} \leq YY_{jktl}M  \quad\qquad \qquad \forall  j,k,t,l,s  \label{eq9: }\\
    Z_{kmtl}^{s} \leq ZZ_{kmtl}M  \qquad \qquad \forall k,m,t,l,s \label{eq10: }
\end{align}

Constraints \eqref{eq11: }-\eqref{eq13: } guarantee that one transportation mode is used to ship a batch of product between different echelons.
\begin{align}
    \sum_{l\in L}XX_{ijtl} \leq 1  \qquad \qquad \forall i,j,t \label{eq11: }\\
    \sum_{l\in L}YY_{jktl} \leq 1  \qquad \qquad \forall j,k,t \label{eq12: }\\
    \sum_{l\in L}ZZ_{jktl} \leq 1  \qquad \qquad \forall k,m,t \label{eq13: }
\end{align}

Constraints \eqref{eq14: }, \eqref{eq15: }, and \eqref{eq16: } show the inventory balance of manufacturers, warehouses, and retailers, respectively. 
\begin{align}
   \sum_{i,i^{'}\in I}\sum_{l\in L}X_{ijtl}^{s}+MI_{jt-1}^{s}+MS_{jt}^{s} =  \sum_{k,k^{'}\in K}\sum_{l\in L}Y_{jktl}^{s}+MI_{jt}^{s}+MS_{jt-1}^{s} \qquad  \forall j,t,s \quad\label{eq14: }\\
   \sum_{j,j^{'}\in J}\sum_{l\in L}Y_{jktl}^{s} =  \sum_{m,m^{'}\in M}\sum_{l\in L}Z_{kmtl}^{s} \qquad \qquad \forall k,t,s \qquad\qquad\qquad\qquad\qquad\quad \quad\quad\quad\label{eq15: }\\
   \sum_{k,k^{'}\in K}\sum_{l\in L}Z_{kmtl}^{s} =  D_{mt}^{s} \qquad \qquad \forall m,t,s \qquad\qquad\qquad\qquad\qquad\qquad\qquad\qquad\quad \quad\label{eq16: }
\end{align}

Constraint \eqref{eq17: } enforces the defined $CO2$ emission cap to the manufacturers.
\begin{equation}
\begin{split}   
  \sum_{j,j^{'}\in J} \sum_{k,k^{'}\in K} \sum_{l\in L} \sum_{s\in S} Y_{jktl}^{s}Cep_{j} + \sum_{i,i^{'}\in I} \sum_{j,j^{'}\in J} \sum_{l\in L} \sum_{s\in S} X_{ijtl}^{s}Cet1_{ijl}\qquad\qquad\qquad\qquad\qquad\\ + \sum_{j,j^{'}\in J} \sum_{k,k^{'}\in K} \sum_{l\in L} \sum_{s\in S} Y_{jktl}^{s}Cet2_{jkl} + \sum_{k,k^{'}\in K} \sum_{m,m^{'}\in M} \sum_{l\in L} \sum_{s\in S} Z_{kmtl}^{s}Cet3_{kml} \leq Cap_{t} \qquad  \forall t \label{eq17: }
\end{split}
\end{equation}

Constraint \eqref{eq18: } describes the minimum number of suppliers that must be selected when multiple sourcing strategy is applied.
\begin{equation}
\begin{split}   
  \sum_{i,i^{'}\in J}XX_{ijtl} \geq Sm  \qquad \qquad \forall j,t,l \label{eq18: }
\end{split}
\end{equation}

Constraints \eqref{eq19: } and \eqref{eq20: } calculate the capacity of manufacturers and warehouses under normal condition as well as under disruption.
\begin{equation}
\begin{split}   
  \sum_{j,j^{'}\in J}Y_{jktl}^{s} \leq \sum_{j,j^{'}\in J}Mc_{j}(1-r_{j}^{s})  \qquad \qquad \forall j,t,l \label{eq19: }
\end{split}
\end{equation}
\begin{equation}
\begin{split}   
  \sum_{k,k^{'}\in K}Z_{kmtl}^{s} \leq \sum_{k,k^{'}\in K}Wc_{k}(1-r_{k}^{'s})  \qquad \qquad \forall j,t,l \label{eq20: }
\end{split}
\end{equation}

\subsection{Resilient strategies}\label{Model:Resilient strategies}

\subsubsection{Backup suppliers}

The backup supplier strategy examines the case of hiring extra possible suppliers. Although hiring extra suppliers imposes higher costs to the supply chain, these suppliers will be available with a higher probability in case of a disruption. By considering this strategy, the range of index $i$ expands to $\{1,2,\cdots,I,I+1,\cdots,I^{'}\}$, where regular suppliers are represented by $\{1,2,\cdots,I\}$, and suppliers $\{I+1,I+2,\cdots,I^{'}\}$ portray the backup suppliers. In the event of a disruption, which may cause delay in suppliers' orders or reduction of their capacity, backup suppliers may curtail interruption.

\subsubsection{Multiple sourcing}

Multiple sourcing diversifies the suppliers, so the supply chain does not heavily rely on a limited number of suppliers, as they may be unavailable when a disruption occurs. Adding constraint \eqref{eq21: } to the model conduces this resilient strategy to our RGSCD problem. Equation \eqref{eq21: } ensures that a minimum number of suppliers are selected to procure the required raw material.
\begin{equation}
\begin{split}   
  \sum_{i,i^{'}\in I}XX_{ijtl} \geq Sm  \qquad \qquad \forall j,t,l \label{eq21: }
\end{split}
\end{equation}

\subsubsection{Safety stock}

Bearing safety stock is another risk management strategy that ensures capability of the supply chain in immediately responding to customer requests, even if delays happen while transporting raw material or final products between different supply chain echelons. Keeping safety stock is accomplished via two adjustments to the original model. First, a maintenance cost for the safety stock equal to $\sum_{j,j^{'}\in J}\sum_{t \in T}\sum_{s\in S}MSS_{jt}^{s}Mic_{j}^{s}$ is added to the inventory costs. Furthermore, the extra inventory held as safety stock in form of final product in each time period is added to equation \eqref{eq14: } to update the manufacturers' inventory balance. After applying this strategic decision to the model, equation \ref{eq14: } is updated as follows.
\begin{equation}
\begin{split}   
  \sum_{i,i^{'}\in I}\sum_{l\in L}X_{ijtl}^{s}+MI_{jt-1}^{s}+MS_{jt}^{s} +MSS_{jt-1}^{s}=  \sum_{k,k^{'}\in K}\sum_{l\in L}Y_{jktl}^{s}+MI_{jt}^{s}\\+MS_{jt-1}^{s}+MSS_{jt}^{s} \qquad \qquad \qquad\qquad\forall j,t,s \qquad\qquad\qquad\qquad\qquad\label{eq22: }
\end{split}
\end{equation}

\subsubsection{Stockpiling}

Stockpiling is an undertaking in which, during the disruption, the manufacturers respond to a proportion of demand using the inventory which had been produced and held before the disruption occurred. This strategy lowers the chance of shortage and is extremely valuable, but incurs additional costs including considerably higher holding costs. Generally, stockpiling does not provide measurable benefits prior to disruptions \citep{liu2016building}. Using the stockpiling system adds a cost equal to $\sum_{j\in J}\sum_{t\in T}\sum_{s\in S}SP_{jt}^{s}Spp$ to the cost function. Also, the manufacturers' inventory balance constraint, equation \eqref{eq14: }, changes as follows.
\begin{equation}
\begin{split}   
  \sum_{i,i^{'}\in I}\sum_{l\in L}X_{ijtl}^{s}+MI_{jt-1}^{s}+MS_{jt}^{s} +MSS_{jt-1}^{s} + SP_{jt}^{s}=  \sum_{k,k^{'}\in K}\sum_{l\in L}Y_{jktl}^{s}+MI_{jt}^{s}\\+MS_{jt-1}^{s}+MSS_{jt}^{s} \qquad \qquad\qquad \qquad \forall j,t,s \qquad \qquad\qquad \qquad\qquad \qquad\quad\quad \label{eq23: }
\end{split}
\end{equation}

\subsubsection{Temporary facilities}

Temporary facilities expedite the supply chain's recovery from a structural disruption \citep{yilmaz2021ensuring}. When needed, possible temporary facilities including manufacturing centers and warehouses act as emergency response centers. To employ this option, domains of indices $i$ and $j$ expand to $\{1, \cdots, I, I+1, \cdots, I^{'}\}$ and $\{1, \cdots, J, J+1, \cdots, J^{'}\}$, respectively. This strategy is useful when the capacity of facilities is diminished, which is commonly observed during disruption events.

\subsubsection{Information sharing system}

One of the main consequences of disruptions is delay. Using an efficient system to share the information between all supply chain members is crucial. For instance, a manufacturer can monitor the demand received by retailers to check the possible spikes, and avoid delays by increasing its production capacity and providing more vehicles to transport the final products. Even though using this information sharing system is costly, it helps the supply chain decrease the cost of delay. Applying this system in our model adds $XX_{ijtl}(So+Tr)$ to the cost function, but reduces transportation delay between different echelons in the amounts specified by $St1_{ijl}$, $St2_{jkl}$, and $St3_{kml}$.

\section{Numerical experiment}\label{Numerical experiment}

\textbf{Data}

To investigate the impact of disruption on supply chain and find the best strategic decisions to control the ripple effect, a numerical experiment was performed using Python programming language and the data extracted from \cite{ghomi2018fuzzy}. Note that since some parameters of the proposed model do not exist in the models presented by \cite{ghomi2018fuzzy}, their values were simulated as will be discussed below. This \attachfile{Dataset_final.xlsx} table (double click on the link to open the file) summarizes the data. Four scenarios were considered to represent different levels of disruption in the model. The first scenario belongs to the normal situation; the rest of the scenarios epitomize three disruption levels, from weak to strong. Three backup suppliers, as well as three temporary facilities for the manufacturing stage and warehouses were considered to be used in the event of a disruption, which is destined to increase the transportation costs, delays, inventory and shortage costs, demand for essential products, and decrease the capacity of manufacturers and warehouses. It is assumed that the government gradually decreases the emission cap through time.

\textbf{Safety stock VS. stockpiling}

\begin{figure}[h]
\centering
\includegraphics[width=10cm]{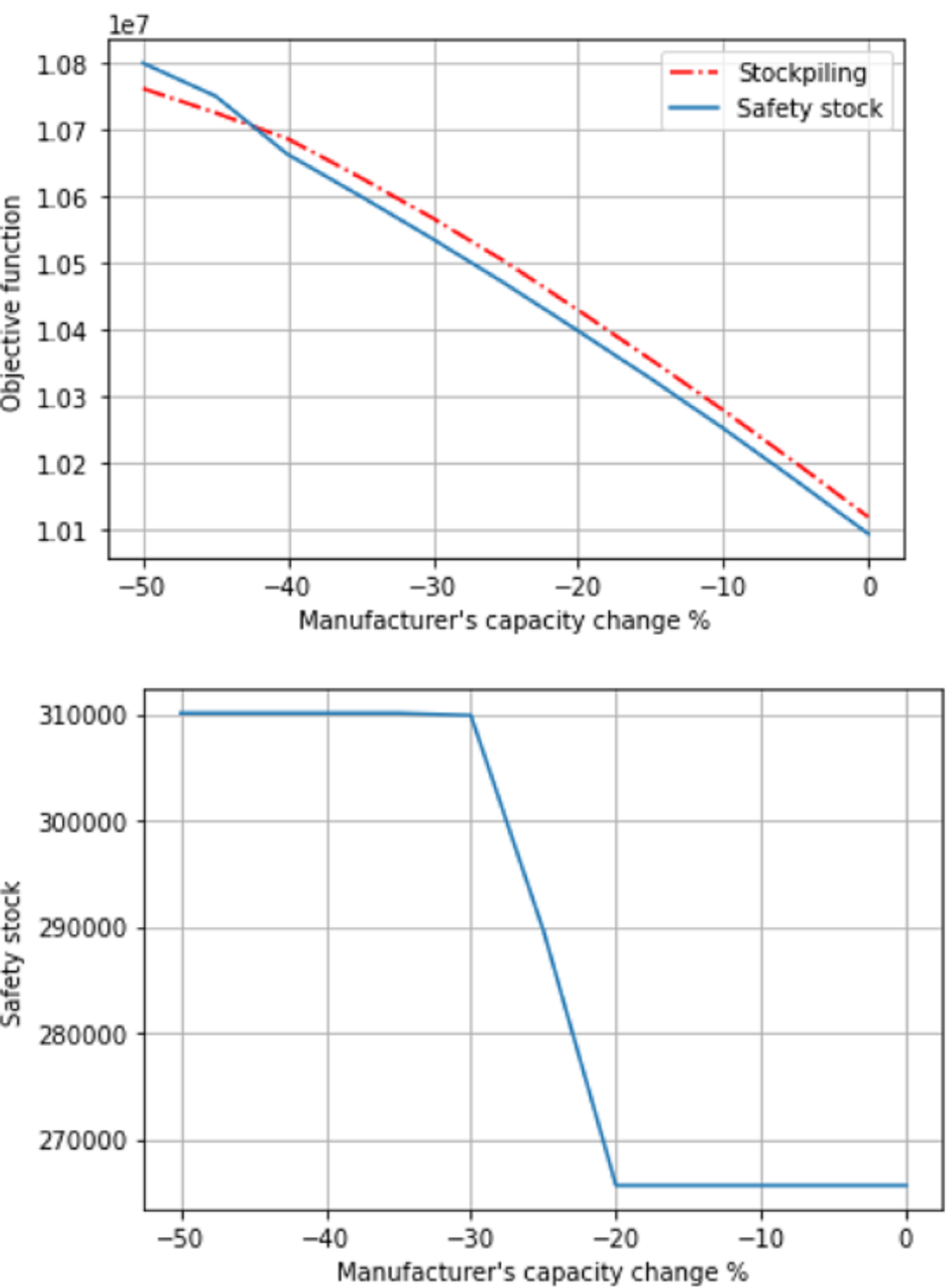}
\caption{Stockpiling vs. safety stockpiling}
\label{Fig: Stockpiling vs. safety stock}
\end{figure}

Stockpiling and safety stock were both used to handle the ripple effect as well as demand spikes. Although from different sources, these two strategies are considered to be identical in controlling disruption effects by procuring the needed products. Thus, they are compared so the managers can choose the best strategy to minimize the costs. Using the presented data set, the developed model was solved considering either the stockpiling or the safety stock strategies under different levels of manufacturers' capacity disruption. The results show that safety stock is a better option when the manufacturers' production capacity drops below 44\%. Figure \ref{Fig: Stockpiling vs. safety stock} depicts the conditions under which each one of the two mentioned strategies are dominant. As it is shown in figure \ref{Fig: Stockpiling vs. safety stock}, by decreasing the manufacturers' capacity, the optimal amount of safety stock the manufacturers should keep increases, which is helpful for decision-makers in their pre-disruption strategy preparations.

\textbf{Backup supplier VS. Multiple sourcing}

\begin{figure}[h]
\centering
\includegraphics[width=10cm]{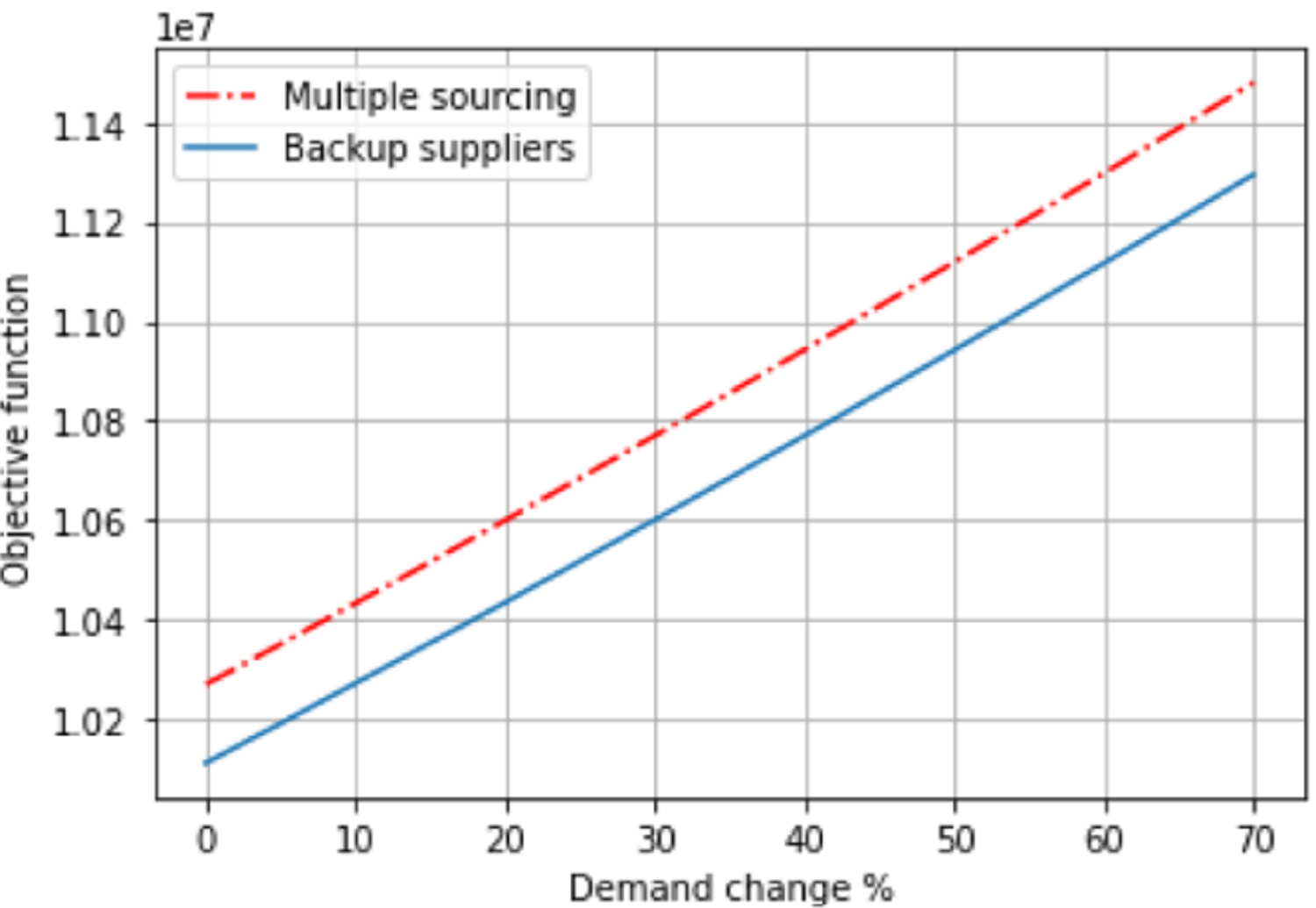}
\caption{Multiple sourcing vs. backup supplier}
\label{Fig: Multiple sourcing vs. backup supplier}
\end{figure}

Two other beneficial pre-disruption strategies in controlling the ripple effect are using backup suppliers and multiple sourcing. Backup suppliers, though costly, prepare the supply chain for disruption in the supply side, and help with procuring raw material and component parts by lowering the level of risk. Multiple sourcing strategy achieves the same results by diversifying the suppliers and avoiding reliance on a limited group of suppliers. The multiple sourcing strategy decreases the supply chain risk, but it has certain disadvantages such as higher costs, lower material quality, and longer lead time due to involvement of secondary suppliers alongside the best selected suppliers in procuring the required material for the manufacturers. In this section, these two strategies are compared to find the optimal and most resilient strategy. Based on figure \ref{Fig: Multiple sourcing vs. backup supplier}, using backup suppliers outperforms multiple sourcing in achieving a higher total utility.

\textbf{Results of the final model}

After comparing similar strategies and determining the dominant game plans using the presented data set, a comprehensive RGSCD model is solved based on four resilient strategies: a) safety stock; b) backup suppliers; c) temporary facilities; and, d) information sharing system. As presented in table \ref{Table: Results}, the best combination of supply chain partners, optimal product flow between different supply chain echelons, and desired amount of manufacturers' safety stock, inventory, and shortage in each period are determined. $Z_{Total}$ shown in table \ref{Table: Results} is the best total utility obtained under all scenarios. This objective is derived from a combination of strategies that help supply chains recover themselves from disruptions. Manufacturers are allowed to produce extra products as long as they have enough storage space. Also, they are allowed to have shortage and respond to that shortage in the future periods. In this regard, the optimum amount of manufacturer's storage and shortage at the end of each period are achieved. For instance, for the presented test problem, $MSS_{31}^{1}$ indicates that under scenario 1, manufacturer $3$ should hold $73,874$ safety stocks to be well prepared for demand fluctuations. Furthermore, the table determines backup suppliers that should be used and temporary facilities that should be opened to maintain the service level during disruptions. The last three rows of the table specify the optimal product flow between suppliers and manufacturers, manufacturers and warehouses, and warehouses and retailers.

\begin{landscape}
\begin{table}
    \centering
    \caption{Optimal results of the proposed model}
    \begin{adjustbox}{max width=20cm}
    \renewcommand{\arraystretch}{1}
\newcommand*{\TitleParbox}[1]{\parbox[c]{10.75cm}{\raggedright #1}}%
    \begin{tabular}{|l|l|}
    \hline
        \textbf{Decision variable} & \textbf{Optimal Value} \\ \hline
         Objective functions & $Z_{Total} = 8,005,492.07$\qquad \qquad $Z_{1} = 4,451,676.10$ \qquad \qquad $Z_{2} = 3,553,815.97$ \\ \hline
        Manufacturers' inventory & $MI_{32}^{1}=36,822$ \\ \hline
        Manufacturers' shortage & $MS_{14}^{1}=14,807$\qquad $MS_{15}^{1}=29,665$\qquad $MS_{16}^{1}=44,465$\qquad $MS_{34}^{1}=36,857$\qquad $MS_{35}^{1}=73,779$\qquad 
        $MS_{36}^{1}=110,633$  \\ \hline
        Safety stock & $MSS_{111}=29,604$ \qquad $MSS_{112}=14,780$ \qquad $MSS_{131}=73,874$ \\ \hline
        Selected suppliers and facilities & $XX_{4321}=1$, $XX_{4331}=1$, $XX_{4341}=1$, $XX_{4351}=1$, $XX_{4361}=1$, $XX_{5121}=1$, $XX_{5131}=1$, $XX_{5141}=1$, $XX_{5151}=1$, $XX_{5161}=1$ \\ 
        & $YY_{1521}=1$, $YY_{1531}=1$,$YY_{1541}=1$,$YY_{1551}=1$,$YY_{1561}=1$,$YY_{3421}=1$,$YY_{3431}=1$,$YY_{3441}=1$,$YY_{3451}=1$,$YY_{3461}=1$,$YY_{3721}=1$, $YY_{3731}=1$,$YY_{3741}=1$,$YY_{3751}=1$,$YY_{3761}=1$,$YY_{3821}=1$,$YY_{3831}=1$,$YY_{3841}=1$,\\& $YY_{3851}=1$,$YY_{3861}=1$  \\ & 
         $ ZZ_{4321}=1$, $ZZ_{4331}=1$, $ZZ_{4341}=1$, $ZZ_{4351}=1$, $ZZ_{4361}=1$, $ZZ_{4721}=1$, $ZZ_{4731}=1$,$ ZZ_{4741}=1$, $ZZ_{4751}=1$, $ZZ_{4761}=1$, $ZZ_{5421}=1$, $ZZ_{5431}=1$, $ZZ_{5441}=1$, $ZZ_{5451}=1$, $ZZ_{5461}=1$, $ZZ_{5621}=1$, $ZZ_{5631}=1$, $ZZ_{5641}=1$,\\ & $ZZ_{5651}=1$, $ZZ_{5661}=1$, $ZZ_{7121}=1$, $ZZ_{7131}=1$, $ZZ_{7141}=1$, $ZZ_{7151}=1$, $ZZ_{7161}=1$,$ ZZ_{8221}=1$, $ZZ_{8231}=1$, $ZZ_{8241}=1$, $ZZ_{8251}=1$, $ZZ_{8261}=1$, $ZZ_{9121}=1$,$ ZZ_{9131}=1$, $ZZ_{9141}=1$, $ZZ_{9151}=1$, $ZZ_{9161}=1$, $ZZ_{9421}=1$,\\ & $ZZ_{9431}=1$, $ZZ_{9441}=1$, $ZZ_{9451}=1$, $ZZ_{9461}=1$, $ZZ_{9521}=1$, $ZZ_{9531}=1$, $ZZ_{9541}=1$, $ZZ_{9551}=1$, $ZZ_{9561}=1$,$ ZZ_{10721}=1$, $ZZ_{10731}=1$, $ZZ_{10741}=1$, $ZZ_{10751}=1$, $ZZ_{10761}=1$  \\ \hline
        Product flow from suppliers to manufacturers & $X_{4321}^{2}=38,643$, $X_{24331}^{2}=38,643$, $X_{4341}^{2}=39,189$, $X_{4351}^{2}=38025$, $X_{4361}^{2}=38,632$, $X_{5121}^{2}=15,671$, $X_{5131}^{2}=15,477$, $X_{5141}^{2}=15,682$, $X_{5151}^{2}=15,720$, $X_{5161}^{2}=15,822$, $X_{4321}^{3}=41,134$, $X_{4331}^{3}=41,717$, $X_{4341}^{3}=41,464$,\\ & $X_{4351}^{3}=41,259$, $X_{4361}^{3}=41,127$,  
        $X_{5121}^{3}=16,578$, $X_{5131}^{3}=16,439$, $X_{5141}^{3}=16,428$, $X_{5151}^{3}=16,909$, $X_{5161}^{3}=16,495$, $X_{4321}^{4}=48,711$, $X_{4331}^{4}=46,704$, $X_{4341}^{4}=48,284$, $X_{4351}^{4}=46,340$, $X_{4361}^{4}=47,417$,  $X_{5121}^{4}=19,510$,\\& $X_{5131}^{4}=18,789$, $X_{5141}^{4}=18,27$, $X_{5151}^{4}=18,550$, $X_{5161}^{4}= 18,305$  \\ \hline
        Product flow from manufacturers to warehouses & $Y_{1521}^{1}=14,824$, $Y_{1531}^{1}=14,780$, $Y_{1541}^{1}=14,807$, $Y_{1551}^{1}=14,858$, $Y_{1561}^{1}=14,800$, 
        $Y_{3421}^{1}=14,871$, $Y_{3431}^{1}=14,743$, $Y_{3441}^{1}=14,715$,  $Y_{3451}^{1}=14,794$, $Y_{3461}^{1}=14,825$,  
        $Y_{3721}^{1}=7,311$, $Y_{3731}^{1}=7,405$, $Y_{3741}^{1}=7,336$, $Y_{3751}^{1}=7,395$, \\& $Y_{3761}^{1}=7.310$,
        $Y_{3821}^{1}=14,870$,  $Y_{3831}^{1}=14,674$, $Y_{3841}^{1}=14,806$, $Y_{3851}^{1}=14,733$, $Y_{3861}^{1}=14,719$, 
        $Y_{1521}^{2}=15,671$, $Y_{1531}^{2}=15,477$, $Y_{1541}^{2}=15,682$, $Y_{1551}^{2}=15,720$,$Y_{1561}^{2}=15,822$,  
        $Y_{3421}^{2}=15,721$, $Y_{3431}^{2}=15,599$, $Y_{3441}^{2}=15,825$, \\ & $Y_{3451}^{2}=15,218$, $Y_{3461}^{2}=15,216$,  
        $Y_{3721}^{2}=7,500$, $Y_{3731}^{2}=7,835$, $Y_{3741}^{2}=7,795$, $Y_{3751}^{2}=7,639$, $Y_{3761}^{2}=7,893$,  
        $Y_{3821}^{2}=15,486$, $Y_{3831}^{2}=15,209$, $Y_{3841}^{2}=15,569$, $Y_{3851}^{2}=15,168$, $Y_{3861}^{2}=15,523$,  
        $Y_{1521}^{3}=16,578$, $Y_{1531}^{3}=16,439$,\\ & $Y_{1541}^{3}=16,428$, $Y_{1551}^{3}=16,909$, $Y_{1561}^{3}=16,495$,  
        $Y_{3421}^{3}=16,550$, $Y_{3431}^{3}=16,523$, $Y_{3441}^{3}=16,716$,  $Y_{3451}^{3}=16,718$, $Y_{3461}^{3}=16,492$,  
        $Y_{3721}^{3}=8,736$, $Y_{3731}^{3}=8,357$, $Y_{3741}^{3}=8,220$, $Y_{3751}^{3}=8,340$, $Y_{3761}^{3}=8,298$,  
        $Y_{3821}^{3}=16,388$,\\ & $Y_{3831}^{3}=16,837$, $Y_{3841}^{3}=16,528$, $Y_{3851}^{3}=16,201$, $Y_{3861}^{3}=16,337$,  
        $Y_{1521}^{4}=19,510$, $Y_{1531}^{4}=18,789$, $Y_{1541}^{4}=18,227$, $Y_{1551}^{4}=18,550$, $Y_{1561}^{4}=18,305$,  
        $Y_{3421}^{4}=19,509$, $Y_{3431}^{4}=18,802$, $Y_{3441}^{4}=18,969$, $Y_{3451}^{4}=16,918$, $Y_{3461}^{4}=19,166$  \\ &
        $Y_{3721}^{4}=10,234$, $Y_{3731}^{4}=9,111$, $Y_{3741}^{4}=10,121$, $Y_{3751}^{4}=10,303$, $Y_{3761}^{4}=10,225$,  
        $Y_{3821}^{4}=18,968$, $Y_{3831}^{4}=18,791$, $Y_{3841}^{4}=19,194$, $Y_{3851}^{4}=19,119$, $Y_{3861}^{4}=18,026$  \\ \hline
        Product flow from warehouses to retailers &  $Z_{4321}^{1}=7,439$, $Z_{4331}^{1}=7,317$, $Z_{4341}^{1}=7,378$, $Z_{4351}^{1}=7,311$, $Z_{4361}^{1}=7,417$, $Z_{4721}^{1}=7,432$, $Z_{4731}^{1}=7,426$, $Z_{4741}^{1}=7,337$,  $Z_{4751}^{1}=7,483$, $Z_{4761}^{1}=7,408$, $Z_{5421}^{1}=7,492$, $Z_{5431}^{1}=7,421$, $Z_{5441}^{1}=7,364$, $Z_{5451}^{1}=7,378$, $Z_{5461}^{1}=7,453$\\ & $Z_{5621}^{1}=7,332$, $Z_{5631}^{1}=7,359$, $Z_{5641}^{1}=7,443$, $Z_{5651}^{1}=7,480$, $Z_{5661}^{1}=7,347$, $Z_{7121}^{1}=7,311$, $Z_{7131}^{1}=7,405$, $Z_{7141}^{1}=7,336$, $Z_{7151}^{1}=7,395$, $Z_{7161}^{1}=7,310$, $Z_{8221}^{1}=7,411$, $Z_{8231}^{1}=7,318$, \\ & $Z_{8241}^{1}=7,393$, $Z_{8251}^{1}=7,335$, $Z_{8261}^{1}=7,400$, $Z_{8521}^{1}=7,459$, $Z_{8531}^{1}=7,356$, $Z_{8541}^{1}=7,413$, $Z_{8551}^{1}=7,398$, $Z_{8561}^{1}=7,319$, $Z_{4321}^{2}=7,806$, $Z_{4331}^{2}=7,871$, $Z_{4341}^{2}=7,957$, $Z_{4351}^{2}=7,701$, $Z_{4361}^{2}=7,608$, $Z_{4721}^{2}=7,915$, $Z_{4731}^{2}=7,728$, \\ & $Z_{4741}^{2}=7,868$, $Z_{4751}^{2}=7,517$, $Z_{4761}^{2}=7,608$, $Z_{5421}^{2}=7,975$, $Z_{5431}^{2}=7,622$,$ Z_{5441}^{2}=7,740$, $Z_{5451}^{2}=7,983$, $Z_{5461}^{2}=7,907$, $Z_{5621}^{2}=7,696$, $Z_{5631}^{2}=7,855$, $Z_{5641}^{2}=7,942$, $Z_{5651}^{2}=7,737$, $Z_{5661}^{2}=7,915$, $Z_{7121}^{2}=7,500$, $Z_{7131}^{2}=7,835$, \\ & $Z_{7141}^{2}=7,795$, $Z_{7151}^{2}=7,639$, $Z_{7161}^{2}=7,893$, $Z_{8221}^{2}=7,904$, $Z_{8231}^{2}=7,507$, $Z_{8241}^{2}=7,677$, $Z_{8251}^{2}=7,612$, $Z_{8261}^{2}=7,758$, $Z_{8521}^{2}=7,582$, $Z_{8531}^{2}=7,702$, $Z_{8541}^{2}=7,892$, $Z_{8551}^{2}=7,556$, $Z_{8561}^{2}=7,765$, $Z_{4321}^{3}=8,365$, $Z_{4331}^{3}=8,489$, \\ & $Z_{4341}^{3}=8,279$, $Z_{4351}^{3}=8,461$, $Z_{4361}^{3}=8,335$, $Z_{4721}^{3}=8,185$, $Z_{4731}^{3}=8,034$, $Z_{4741}^{3}=8,437$, $Z_{4751}^{3}=8,257$, $Z_{4761}^{3}=8,157$, $Z_{5421}^{3}=8,220$, $Z_{5431}^{3}=8,173$, $Z_{5441}^{3}=8,108$, $Z_{5451}^{3}=8,428$, $Z_{5461}^{3}=8,213$, $Z_{5621}^{3}=8,358$, $Z_{5631}^{3}=8,266$, \\ & $Z_{5641}^{3}=8,320$, $Z_{5651}^{3}=8,481$, $Z_{5661}^{3}=8,282$, $Z_{7121}^{3}=8,376$, $Z_{7131}^{3}=8,357$, $Z_{7141}^{3}=8,220$, $Z_{7151}^{3}=8,340$, $Z_{7161}^{3}=8,298$, $Z_{8221}^{3}=8,096$, $Z_{8231}^{3}=8,380$, $Z_{8241}^{3}=8,377$, $Z_{8251}^{3}=8,133$, $Z_{8261}^{3}=8,041$, $Z_{8521}^{3}=8,292$, $Z_{8531}^{3}=8,457$, \\ &$Z_{8541}^{3}=8,151$, $Z_{8551}^{3}=8,068$, $Z_{8561}^{3}=8,296$, $Z_{4321}^{4}=9,283$, $Z_{4331}^{4}=9,653$, $Z_{4341}^{4}=10,209$, $Z_{4351}^{4}=6,849$, $Z_{4361}^{4}=8,975$, $Z_{4721}^{4}=10,226$, $Z_{4731}^{4}=9,149$, $Z_{4741}^{4}=8,760$, $Z_{4751}^{4}=10,069$, $Z_{4761}^{4}=10,191$ \\ & $Z_{5421}^{4}=10,002$, $Z_{5431}^{4}=9,902$, $Z_{5441}^{4}=8,760$,$Z_{5451}^{4}=8,953$, $Z_{5461}^{4}=9,047$, $Z_{5621}^{4}=9,508$, $Z_{5631}^{4}=8,887$, $Z_{5641}^{4}=9,467$, $Z_{5651}^{4}=9,597$, $Z_{5661}^{4}=9,258$, $Z_{7121}^{4}=10,234$, $Z_{7131}^{4}=9,111$, 
        \\ & $Z_{7141}^{4}=10,121$, $Z_{7151}^{4}=10,303$, $Z_{7161}^{4}=10,225$, $Z_{8221}^{4}=10,300$, $Z_{8231}^{4}=8,898$, $Z_{8241}^{4}=9,872$, $Z_{8251}^{4}=8,748$, $Z_{8261}^{4}=9,256$, $Z_{8521}^{4}=8,668$, $Z_{8531}^{4}=9,893$, $Z_{8541}^{4}=9,322$, $Z_{8551}^{4}=10,371$, $Z_{8561}^{4}=8,770$ \\
         \hline      
         \end{tabular}
         \end{adjustbox}
    \label{Table: Results}
\end{table}
\end{landscape}

To validate the proposed model and solution approach, and to investigate the supply chain behavior, further analyses are conducted as follows.

\textbf{Behaviour of essential and non-essential products}

\begin{figure}[h]
\centering
\includegraphics[width=\textwidth]{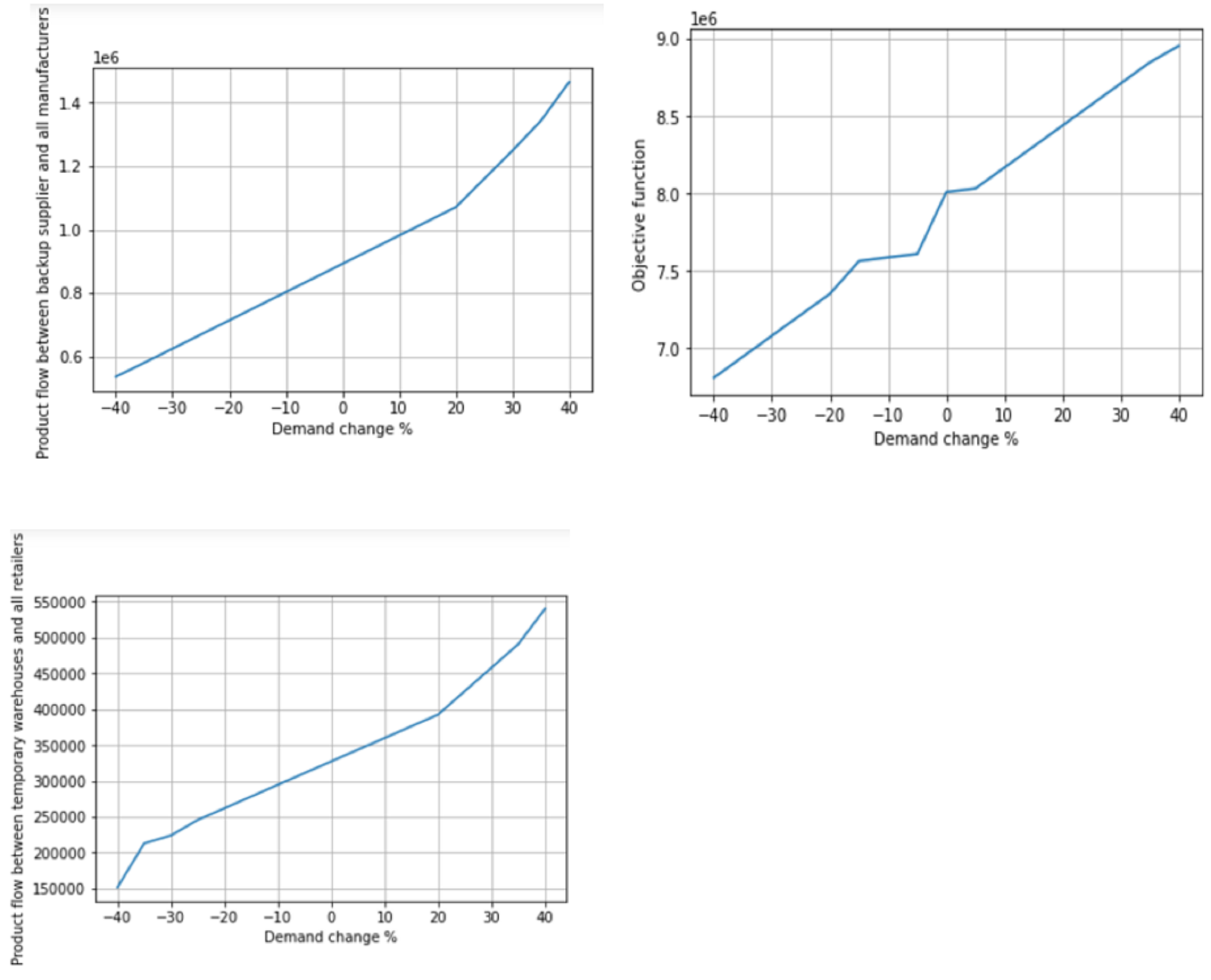}
\caption{Sensitivity analysis on demand fluctuation for essential and non-essential products}
\label{Fig: Demand effect}
\end{figure}

Distinguishing the supply chain's products between essential and non-essential is crucial for the managers. As it was evident during the COVID-19 pandemic, essential products experienced a large spike in demand, which was mostly non-present for non-essential products. To measure the impact of demand change on the optimal solution of the designed model a sensitivity analysis was designed. For this purpose, different rates of demand fluctuation were considered and the results were analyzed. Demand increases can be interpreted as the fluctuation of demand for essential products; demand decrease represents the market uncertainties for non-essential products. Figure \ref{Fig: Demand effect} demonstrates that demand intensification calls for expanding the adoption of temporary facilities for manufacturers and warehouses. Also, based on figure \ref{Fig: Demand effect}, elevation of demand induces higher need for product flow between temporary facilities, which indicates a desired response to disruptions via activation of the temporary facilities.

\textbf{The effects of carbon abatement regulation}

\begin{figure}[h]
\centering
\includegraphics[width=10cm]{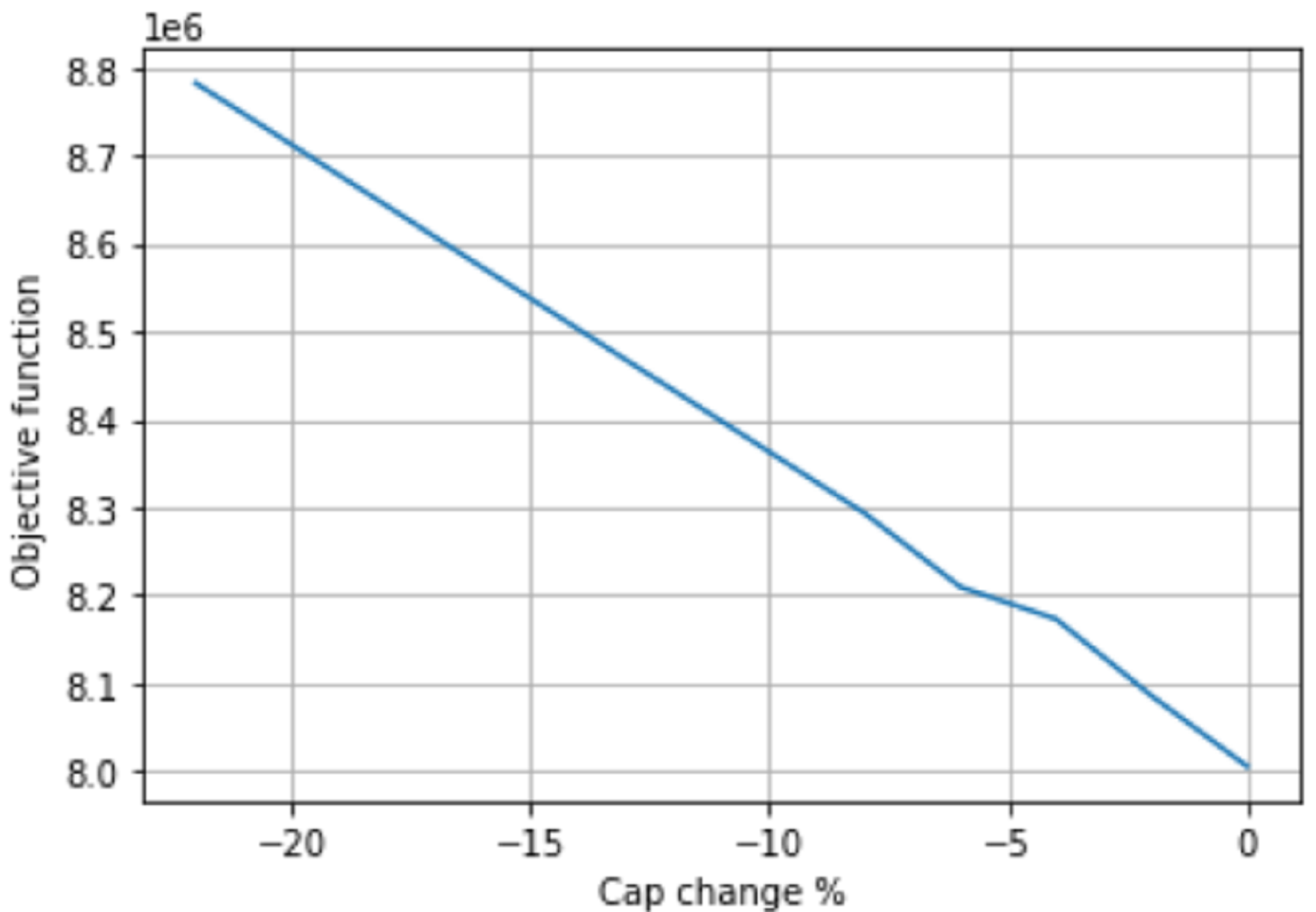}
\caption{The effect of carbon abatement regulations}
\label{Fig: Cap effect}
\end{figure}

An array of carbon emission regulations are employed in various countries to restrict supply chains from harming the environment. One way to restrict the emission levels is to place a cap on the manufacturers' emission. Accordingly, the proposed model assumes that a cap-and-trade regulation is enforced, where the cap is gradually decreased by governments to reduce total $CO_2$ emissions and to allow manufacturers to adapt to the regulation requirements. Hence, we analyze the impact of cap modifications on the supply chain in this section. As it is shown in figure \ref{Fig: Cap effect}, reducing the cap imposes enormous costs to the supply chain for $CO_2$ emissions; the model reacts by reducing the emission level, which results in total utility contraction. After reducing the cap to 76\% of its original amount, the model becomes infeasible. In other words, reducing the cap beyond a certain threshold makes it impossible to respond to customers' demand, while adhering to the maximum emission requirements. Therefore, the government is advised to recognize supply chain limitations before setting the cap amounts. The other methods of reducing the emission levels without decreasing the production level is discussed in \cite{mirzaee2022three}.

\textbf{Capacity decrease effect}

\begin{figure}[h]
\centering
\includegraphics[width=10cm]{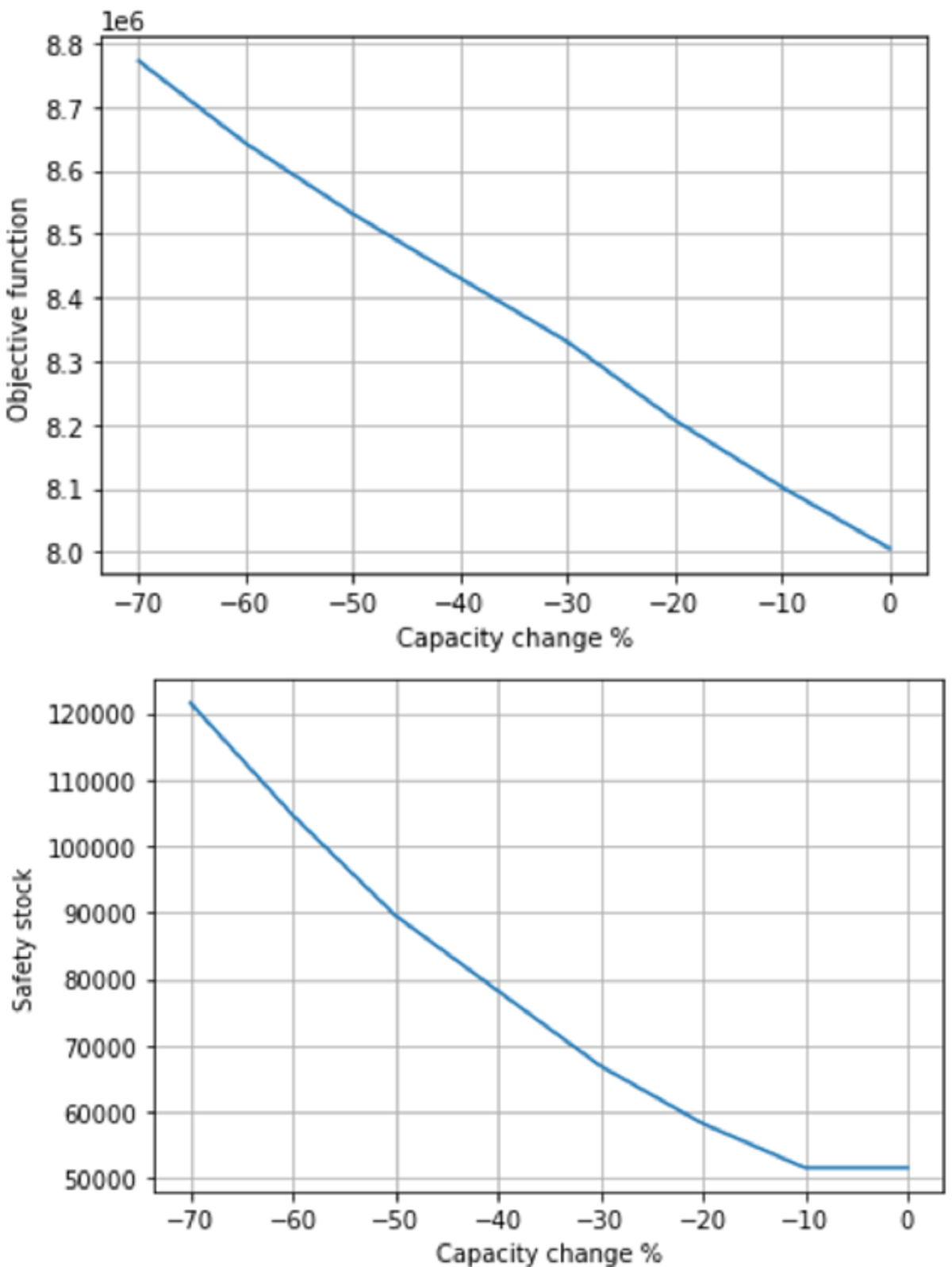}
\caption{Capacity effect}
\label{Fig: Capacity effect}
\end{figure}

Ripple effect usually impacts the capacities of manufacturers and warehouses \citep{monostori2021mitigation}, which calls for an appropriate risk management strategy, such as responding to customers' demand by keeping more safety stock. A sensitivity analysis is conducted to assess the correct strategies when the capacities are decreased. As figure \ref{Fig: Capacity effect} illustrates, decreasing the capacity of the supply chain facilities deteriorates the objective function value for it forces the model to responds to demand by keeping more safety stock, which is more expensive to carry. Decreasing the capacity beyond a certain threshold causes infeasibility as the limited capacity makes it impossible to handle the realized demand. Consequently, the supply chain managers must consider keeping an appropriate amount of temporary facilities available to prevent service level degradation.

\section{Conclusion}\label{Conclusion}

This paper proposes a resilient green supply chain design to mitigate the ripple effects of disruptions such as those experienced during the COVID-19 pandemic, and contemplates the best design strategies to alleviate the impacts of turmoil. Accordingly, six resilience strategies were applied to the designed supply chain to prevent service level reduction. Also, a two-stage stochastic optimization approach was employed to control the ripple effect and the inherent uncertainty in parameter estimation. In other words, this study provides a decision making framework for supply chain managers to use the best transportation channel for materials and final products, enhance service level, and control uncertainty of estimating model parameters and disruptions.

The validity of the designed model and proposed solution approach were verified through a numerical example. The results confirmed that the selected strategies mitigate the ripple effect in presence of stochastic parameters and possible disruptions. The analysis of demand fluctuations and facility capacity variations allow the supply chain managers make the best possible strategic decisions for their procurement and production plans before facing disruptions. Moreover, sensitivity analysis of carbon cap reveals that governments must provide a proper amount of carbon allowance for manufacturers.

This paper pioneers the study of resilient green supply chain design using two-stage stochastic optimization approach considering a variety of pre- and post-disruption resilience strategies under carbon emission regulations. Given the novelty of this issue and its importance after the COVID-19 pandemic, additional studies in this field are required. It is worth investigating perishable products with a limited lifetime in this context. Keeping perishable products as inventory or safety stock will be more expensive for the supply chain, and stockpiling or outsourcing may be a better option. Additionally, considering back-to-normal plans to decrease the extra costs that supply chain partners pay during the disruption events is an interesting topic. Furthermore, scrutinizing the relationship between safety stock and service level is another venue worth further investigation. 

\bibliographystyle{apalike}
\bibliography{References.bib}

\begin{thebibliography}{}

\bibitem[Amiri-Aref et~al., 2018]{amiri2018multi}
Amiri-Aref, M., Klibi, W., and Babai, M.~Z. (2018).
\newblock The multi-sourcing location inventory problem with stochastic demand.
\newblock {\em European Journal of Operational Research}, 266(1):72--87.

\bibitem[Badri et~al., 2017]{badri2017two}
Badri, H., Ghomi, S.~F., and Hejazi, T.-H. (2017).
\newblock A two-stage stochastic programming approach for value-based
  closed-loop supply chain network design.
\newblock {\em Transportation Research Part E: Logistics and Transportation
  Review}, 105:1--17.

\bibitem[Bhatia and Gangwani, 2021]{bhatia2021green}
Bhatia, M.~S. and Gangwani, K.~K. (2021).
\newblock Green supply chain management: Scientometric review and analysis of
  empirical research.
\newblock {\em Journal of cleaner production}, 284:124722.

\bibitem[Boskabadi et~al., 2022]{boskabadi2022design}
Boskabadi, A., Mirmozaffari, M., Yazdani, R., and Farahani, A. (2022).
\newblock Design of a distribution network in a multi-product, multi-period
  green supply chain system under demand uncertainty.
\newblock {\em Sustainable Operations and Computers}, 3:226--237.

\bibitem[Cui et~al., 2020]{cui2020stochastic}
Cui, L., Deng, J., Liu, R., Xu, D., Zhang, Y., and Xu, M. (2020).
\newblock A stochastic multi-item replenishment and delivery problem with
  lead-time reduction initiatives and the solving methodologies.
\newblock {\em Applied mathematics and computation}, 374:125055.

\bibitem[Darestani and Hemmati, 2019]{darestani2019robust}
Darestani, S.~A. and Hemmati, M. (2019).
\newblock Robust optimization of a bi-objective closed-loop supply chain
  network for perishable goods considering queue system.
\newblock {\em Computers \& industrial engineering}, 136:277--292.

\bibitem[Fattahi et~al., 2017]{fattahi2017responsive}
Fattahi, M., Govindan, K., and Keyvanshokooh, E. (2017).
\newblock Responsive and resilient supply chain network design under
  operational and disruption risks with delivery lead-time sensitive customers.
\newblock {\em Transportation research part E: Logistics and transportation
  review}, 101:176--200.

\bibitem[Foroozesh et~al., 2022]{foroozesh2022green}
Foroozesh, N., Karimi, B., and Mousavi, S. (2022).
\newblock Green-resilient supply chain network design for perishable products
  considering route risk and horizontal collaboration under robust
  interval-valued type-2 fuzzy uncertainty: A case study in food industry.
\newblock {\em Journal of Environmental Management}, 307:114470.

\bibitem[Gholami-Zanjani et~al., 2021]{gholami2021design}
Gholami-Zanjani, S.~M., Klibi, W., Jabalameli, M.~S., and Pishvaee, M.~S.
  (2021).
\newblock The design of resilient food supply chain networks prone to epidemic
  disruptions.
\newblock {\em International Journal of Production Economics}, 233:108001.

\bibitem[Ghomi-Avili et~al., 2018]{ghomi2018fuzzy}
Ghomi-Avili, M., Naeini, S. G.~J., Tavakkoli-Moghaddam, R., and Jabbarzadeh, A.
  (2018).
\newblock A fuzzy pricing model for a green competitive closed-loop supply
  chain network design in the presence of disruptions.
\newblock {\em Journal of Cleaner Production}, 188:425--442.

\bibitem[Hasani et~al., 2021]{hasani2021multi}
Hasani, A., Mokhtari, H., and Fattahi, M. (2021).
\newblock A multi-objective optimization approach for green and resilient
  supply chain network design: a real-life case study.
\newblock {\em Journal of Cleaner Production}, 278:123199.

\bibitem[Hendricks and Singhal, 2005]{hendricks2005empirical}
Hendricks, K.~B. and Singhal, V.~R. (2005).
\newblock An empirical analysis of the effect of supply chain disruptions on
  long-run stock price performance and equity risk of the firm.
\newblock {\em Production and Operations management}, 14(1):35--52.

\bibitem[Hosseini and Ivanov, 2019]{hosseini2019new}
Hosseini, S. and Ivanov, D. (2019).
\newblock A new resilience measure for supply networks with the ripple effect
  considerations: A bayesian network approach.
\newblock {\em Annals of Operations Research}, pages 1--27.

\bibitem[Hosseini et~al., 2019]{hosseini2019review}
Hosseini, S., Ivanov, D., and Dolgui, A. (2019).
\newblock Review of quantitative methods for supply chain resilience analysis.
\newblock {\em Transportation Research Part E: Logistics and Transportation
  Review}, 125:285--307.

\bibitem[Hosseini et~al., 2020]{hosseini2020ripple}
Hosseini, S., Ivanov, D., and Dolgui, A. (2020).
\newblock Ripple effect modelling of supplier disruption: integrated markov
  chain and dynamic bayesian network approach.
\newblock {\em International Journal of Production Research},
  58(11):3284--3303.

\bibitem[Hosseini-Motlagh et~al., 2019]{hosseini2019reverse}
Hosseini-Motlagh, S.-M., Nouri-Harzvili, M., Choi, T.-M., and Ebrahimi, S.
  (2019).
\newblock Reverse supply chain systems optimization with dual channel and
  demand disruptions: Sustainability, csr investment and pricing coordination.
\newblock {\em Information Sciences}, 503:606--634.

\bibitem[Ivanov and Dolgui, 2022]{ivanov2022stress}
Ivanov, D. and Dolgui, A. (2022).
\newblock Stress testing supply chains and creating viable ecosystems.
\newblock {\em Operations Management Research}, 15(1):475--486.

\bibitem[Ivanov et~al., 2015]{ivanov2015supply}
Ivanov, D., Dolgui, A., and Sokolov, B. (2015).
\newblock Supply chain design with disruption considerations: Review of
  research streams on the ripple effect in the supply chain.
\newblock {\em IFAC-PapersOnLine}, 48(3):1700--1707.

\bibitem[Ivanov et~al., 2019]{ivanov2019impact}
Ivanov, D., Dolgui, A., and Sokolov, B. (2019).
\newblock The impact of digital technology and industry 4.0 on the ripple
  effect and supply chain risk analytics.
\newblock {\em International Journal of Production Research}, 57(3):829--846.

\bibitem[Ivanov et~al., 2016]{ivanov2016disruption}
Ivanov, D., Pavlov, A., Dolgui, A., Pavlov, D., and Sokolov, B. (2016).
\newblock Disruption-driven supply chain (re)-planning and performance impact
  assessment with consideration of pro-active and recovery policies.
\newblock {\em Transportation Research Part E: Logistics and Transportation
  Review}, 90:7--24.

\bibitem[Jabbarzadeh et~al., 2018]{jabbarzadeh2018resilient}
Jabbarzadeh, A., Fahimnia, B., and Sabouhi, F. (2018).
\newblock Resilient and sustainable supply chain design: sustainability
  analysis under disruption risks.
\newblock {\em International Journal of Production Research},
  56(17):5945--5968.

\bibitem[John et~al., 2018]{john2018multi}
John, S.~T., Sridharan, R., Kumar, P.~R., and Krishnamoorthy, M. (2018).
\newblock Multi-period reverse logistics network design for used refrigerators.
\newblock {\em Applied Mathematical Modelling}, 54:311--331.

\bibitem[Kamalahmadi and Mellat-Parast, 2016]{kamalahmadi2016developing}
Kamalahmadi, M. and Mellat-Parast, M. (2016).
\newblock Developing a resilient supply chain through supplier flexibility and
  reliability assessment.
\newblock {\em International Journal of Production Research}, 54(1):302--321.

\bibitem[Khaloie et~al., 2020]{khaloie2020coordinated}
Khaloie, H., Abdollahi, A., Shafie-Khah, M., Anvari-Moghaddam, A., Nojavan, S.,
  Siano, P., and Catal{\~a}o, J.~P. (2020).
\newblock Coordinated wind-thermal-energy storage offering strategy in energy
  and spinning reserve markets using a multi-stage model.
\newblock {\em Applied Energy}, 259:114168.

\bibitem[Kim et~al., 2015]{kim2015supply}
Kim, Y., Chen, Y.-S., and Linderman, K. (2015).
\newblock Supply network disruption and resilience: A network structural
  perspective.
\newblock {\em Journal of operations Management}, 33:43--59.

\bibitem[Liao, 2018]{liao2018reverse}
Liao, T.-Y. (2018).
\newblock Reverse logistics network design for product recovery and
  remanufacturing.
\newblock {\em Applied Mathematical Modelling}, 60:145--163.

\bibitem[Liu et~al., 2016]{liu2016building}
Liu, F., Song, J.-S., and Tong, J.~D. (2016).
\newblock Building supply chain resilience through virtual stockpile pooling.
\newblock {\em Production and Operations management}, 25(10):1745--1762.

\bibitem[Liu et~al., 2022]{liu2022optimization}
Liu, M., Liu, Z., Chu, F., Dolgui, A., Chu, C., and Zheng, F. (2022).
\newblock An optimization approach for multi-echelon supply chain viability
  with disruption risk minimization.
\newblock {\em Omega}, page 102683.

\bibitem[Memon et~al., 2015]{memon2015group}
Memon, M.~S., Lee, Y.~H., and Mari, S.~I. (2015).
\newblock Group multi-criteria supplier selection using combined grey systems
  theory and uncertainty theory.
\newblock {\em Expert Systems with Applications}, 42(21):7951--7959.

\bibitem[Mirzaee et~al., 2022]{mirzaee2022three}
Mirzaee, H., Samarghandi, H., and Willoughby, K. (2022).
\newblock A three-player game theory model for carbon cap-and-trade mechanism
  with stochastic parameters.
\newblock {\em Computers \& Industrial Engineering}, page 108285.

\bibitem[Mohammed et~al., 2017]{mohammed2017cost}
Mohammed, A., Wang, Q., and Li, X. (2017).
\newblock A cost-effective decision-making algorithm for an rfid-enabled hmsc
  network design: A multi-objective approach.
\newblock {\em Industrial Management \& Data Systems}.

\bibitem[Mohebalizadehgashti et~al., 2020]{mohebalizadehgashti2020designing}
Mohebalizadehgashti, F., Zolfagharinia, H., and Amin, S.~H. (2020).
\newblock Designing a green meat supply chain network: A multi-objective
  approach.
\newblock {\em International Journal of Production Economics}, 219:312--327.

\bibitem[Monostori, 2021]{monostori2021mitigation}
Monostori, J. (2021).
\newblock Mitigation of the ripple effect in supply chains: Balancing the
  aspects of robustness, complexity and efficiency.
\newblock {\em CIRP Journal of Manufacturing Science and Technology},
  32:370--381.

\bibitem[Ni et~al., 2018]{ni2018modeling}
Ni, N., Howell, B.~J., and Sharkey, T.~C. (2018).
\newblock Modeling the impact of unmet demand in supply chain resiliency
  planning.
\newblock {\em Omega}, 81:1--16.

\bibitem[Oksuz and Satoglu, 2020]{oksuz2020two}
Oksuz, M.~K. and Satoglu, S.~I. (2020).
\newblock A two-stage stochastic model for location planning of temporary
  medical centers for disaster response.
\newblock {\em International Journal of Disaster Risk Reduction}, 44:101426.

\bibitem[{\"O}z{\c{c}}elik et~al., 2021]{ozccelik2021robust}
{\"O}z{\c{c}}elik, G., Faruk~Y{\i}lmaz, {\"O}., and Bet{\"u}l~Yeni, F. (2021).
\newblock Robust optimisation for ripple effect on reverse supply chain: an
  industrial case study.
\newblock {\em International Journal of Production Research}, 59(1):245--264.

\bibitem[O’Brien, 2013]{o2013fifty}
O’Brien, C. (2013).
\newblock Fifty years of shifting paradigms.
\newblock {\em International Journal of Production Research},
  51(23-24):6740--6745.

\bibitem[Pavlov et~al., 2017]{pavlov2017hybrid}
Pavlov, A., Ivanov, D., Dolgui, A., and Sokolov, B. (2017).
\newblock Hybrid fuzzy-probabilistic approach to supply chain resilience
  assessment.
\newblock {\em IEEE Transactions on Engineering Management}, 65(2):303--315.

\bibitem[Reuters, 2020]{REUTERS2020news}
Reuters (2020).

\bibitem[Rezapour et~al., 2017]{rezapour2017resilient}
Rezapour, S., Farahani, R.~Z., and Pourakbar, M. (2017).
\newblock Resilient supply chain network design under competition: A case
  study.
\newblock {\em European Journal of Operational Research}, 259(3):1017--1035.

\bibitem[Sawik, 2019]{sawik2019disruption}
Sawik, T. (2019).
\newblock Disruption mitigation and recovery in supply chains using portfolio
  approach.
\newblock {\em Omega}, 84:232--248.

\bibitem[Sharma et~al., 2022]{sharma2022managing}
Sharma, M., Joshi, S., Luthra, S., and Kumar, A. (2022).
\newblock Managing disruptions and risks amidst covid-19 outbreaks: role of
  blockchain technology in developing resilient food supply chains.
\newblock {\em Operations Management Research}, 15(1):268--281.

\bibitem[Snyder et~al., 2016]{snyder2016or}
Snyder, L.~V., Atan, Z., Peng, P., Rong, Y., Schmitt, A.~J., and Sinsoysal, B.
  (2016).
\newblock Or/ms models for supply chain disruptions: A review.
\newblock {\em Iie Transactions}, 48(2):89--109.

\bibitem[Tomlin, 2006]{tomlin2006value}
Tomlin, B. (2006).
\newblock On the value of mitigation and contingency strategies for managing
  supply chain disruption risks.
\newblock {\em Management science}, 52(5):639--657.

\bibitem[Tordecilla et~al., 2021]{tordecilla2021simulation}
Tordecilla, R.~D., Juan, A.~A., Montoya-Torres, J.~R., Quintero-Araujo, C.~L.,
  and Panadero, J. (2021).
\newblock Simulation-optimization methods for designing and assessing resilient
  supply chain networks under uncertainty scenarios: A review.
\newblock {\em Simulation modelling practice and theory}, 106:102166.

\bibitem[Tucker et~al., 2020]{tucker2020incentivizing}
Tucker, E.~L., Daskin, M.~S., Sweet, B.~V., and Hopp, W.~J. (2020).
\newblock Incentivizing resilient supply chain design to prevent drug
  shortages: policy analysis using two-and multi-stage stochastic programs.
\newblock {\em IISE Transactions}, 52(4):394--412.

\bibitem[Wang et~al., 2011]{wang2011multi}
Wang, F., Lai, X., and Shi, N. (2011).
\newblock A multi-objective optimization for green supply chain network design.
\newblock {\em Decision support systems}, 51(2):262--269.

\bibitem[Wang et~al., 2021]{wang2021incentive}
Wang, W., Zhang, Y., Zhang, W., Gao, G., and Zhang, H. (2021).
\newblock Incentive mechanisms in a green supply chain under demand
  uncertainty.
\newblock {\em Journal of Cleaner Production}, 279:123636.

\bibitem[Y{\i}lmaz et~al., 2021]{yilmaz2021ensuring}
Y{\i}lmaz, {\"O}.~F., {\"O}z{\c{c}}elik, G., and Yeni, F.~B. (2021).
\newblock Ensuring sustainability in the reverse supply chain in case of the
  ripple effect: A two-stage stochastic optimization model.
\newblock {\em Journal of cleaner production}, 282:124548.

\bibitem[Zahiri et~al., 2018]{zahiri2018design}
Zahiri, B., Jula, P., and Tavakkoli-Moghaddam, R. (2018).
\newblock Design of a pharmaceutical supply chain network under uncertainty
  considering perishability and substitutability of products.
\newblock {\em Information Sciences}, 423:257--283.

\bibitem[Zhang et~al., 2019]{zhang2019fuzzy}
Zhang, S., Zhang, P., and Zhang, M. (2019).
\newblock Fuzzy emergency model and robust emergency strategy of supply chain
  system under random supply disruptions.
\newblock {\em Complexity}, 2019.

\end{thebibliography}

\end{document}